%% file: delamination2.tex
\documentclass[preprint,10pt]{elsarticle}

\usepackage{url}
\usepackage{graphicx}
\usepackage{caption}
\usepackage{subfig}
\usepackage{epsfig}
\usepackage{psfrag}
\usepackage{hyperref}
\usepackage{amsmath,amssymb,amsfonts,amsthm,stmaryrd}
\usepackage{dsfont}
\usepackage{eucal}
\usepackage{mathrsfs}
\usepackage{listings}    
\usepackage{color}
\usepackage{float}
\usepackage{natbib}
\usepackage{booktabs}
\usepackage{tabularx}
\usepackage{setspace}
\usepackage{algorithmicx}
\usepackage{algorithm}
\usepackage{algpseudocode}

\input{pre_commands.tex}

\newcommand{\grbf}[1]{\mbox{\boldmath{$#1$}}}

\floatstyle{ruled}
\newfloat{Program}{thp}{lop}
\floatname{Program}{Program}

\floatstyle{ruled}
\newfloat{Fbox}{thp}{lop}
\floatname{Fbox}{Box}

\newcolumntype{C}{>{\centering\arraybackslash}X}

\definecolor{darkgray}{rgb}{0.95,0.95,0.95}
\definecolor{mygreen}{rgb}{0,0.6,0}

\lstdefinestyle{Matlab}
{
 basicstyle=\footnotesize, numbers=none, numberstyle=\tiny,%
 showstringspaces=false, language=Matlab, escapechar=|,frame=tb,%
 commentstyle=\color{mygreen}
}

\lstdefinestyle{Matlab-num}
{
 basicstyle=\footnotesize, numbers=left, numberstyle=\tiny,%
 showstringspaces=false, language=Matlab, escapechar=|,frame=tb,%
commentstyle=\color{mygreen}
}

\newcommand{\tty}[1]{\textnormal{\texttt{#1}}}
\newcommand{\sym}[1]{\textnormal{\textit{#1}}}

\lstnewenvironment{snippet}[1][]
{
 \lstset{style=Matlab, xleftmargin=5mm, gobble=4, #1}
}
{}

\lstnewenvironment{snippet1}[1][]
{
 \lstset{style=Matlab-num, xleftmargin=5mm, gobble=4, #1}
}
{}

\def\bibsection{\section*{References}}

\begin{document}


\definecolor{MyDarkBlue}{rgb}{1, 0.9, 1}
\lstset{language=Matlab,
       basicstyle=\footnotesize,
       commentstyle=\itshape,
       stringstyle=\ttfamily,
       showstringspaces=false,
       tabsize=2}
\lstdefinestyle{commentstyle}{color=\color{green}}

\theoremstyle{remark}
\newtheorem{thm}{Theorem}[section]
\newtheorem{rmk}[thm]{Remark}


\definecolor{red}{gray}{0}
\definecolor{blue}{gray}{0}


\begin{frontmatter}

\title{Isogeometric cohesive elements for two and three dimensional
composite delamination analysis}       

\author[cardiff]{Vinh Phu Nguyen  \fnref{fn1}}
\author[cardiff]{Pierre Kerfriden  \fnref{fn2}}
\author[cardiff]{St\'{e}phane P.A. Bordas \corref{cor1}\fnref{fn3}}

\cortext[cor1]{Corresponding author}

\address[cardiff]{School of Engineering, Institute of Mechanics and Advanced
Materials, Cardiff University, Queen's Buildings, The Parade, Cardiff \\
CF24 3AA}
\address[weimar]{Institute of Structural Mechanics, Bauhaus-Universit\"{a}t
Weimar, Marienstra\ss{}e 15 99423 Weimar}

\fntext[fn1]{\url nguyenpv@cardiff.ac.uk}
\fntext[fn2]{\url pierre@cardiff.ac.uk}
\fntext[fn3]{\url stephane.bordas@alum.northwestern.edu}

\begin{abstract}
Isogeometric cohesive elements are presented for modeling two and three dimensional 
delaminated composite structures. We exploit the knot insertion algorithm offered by
NURBS (Non Uniform Rational B-splines) to generate cohesive elements along delamination
planes in an automatic fashion. A complete computational framework is presented including 
pre-processing, processing and post-processing. They are explained in details and
implemented in MIGFEM--an open source Matlab Isogemetric Analysis code developed by the authors. 
The composite laminates are modeled using both NURBS solid and shell elements. Several two and three
dimensional examples ranging from standard delamination tests (the mixed mode bending test), 
the L-shaped specimen with a fillet, three dimensional (3D) double cantilever beam and a 3D singly curved
thick-walled laminate are provided.            
To the authors' knowledge, it is the first time that NURBS-based isogeometric analysis for two/three
dimensional delamination modeling is presented. For all examples considered, the
proposed framework outperforms conventional Lagrange finite elements.
\end{abstract}


\end{frontmatter}


\section{Introduction}

Isogeometric analysis (IGA) was proposed by Hughes and his co-workers 
\cite{hughes_isogeometric_2005} in 2005 to reduce the gap between Computer Aided Design (CAD)
and Finite Element Analysis (FEA). The idea is to use CAD technology such B-splines, NURBS (Non Uniform
Rational B-splines), T-splines \etc as basis functions in a finite element (FE) framework. 
Since this seminal paper, a monograph has been published entirely on
the subject  \cite{cottrel_book_2009} and applications have been found
in several fields including structural mechanics, solid
mechanics, fluid mechanics and contact mechanics. 
It should be emphasized that the idea of using CAD technologies in finite elements is not new.
For example in \cite{NME:NME292}, B-splines were used as shape functions in FEM and
subdivision surfaces were adopted to model shells \cite{Cirak_2000}.

Due to the ultra smoothness provided by NURBS basis, IGA has been successfully applied to 
many engineering problems ranging from contact mechanics, see \eg
\cite{temizer_contact_2011,jia_isogeometric_2011,temizer_three-dimensional_2012,
de_lorenzis_large_2011,Matzen201327}, optimisation problems 
\cite{wall_isogeometric_2008,manh_isogeometric_2011,qian_isogeometric_2011,xiaoping_full_2010},
structural mechanics \cite{benson_isogeometric_2010,kiendl_isogeometric_2009,benson_large_2011,
beirao_da_veiga_isogeometric_2012,uhm_tspline_2009,Echter2013170,Benson2013133}, 
structural vibration \cite{cottrell_isogeometric_2006,Hughes20084104,NME:NME4282,Wang2013}, 
to fluids mechanics 
\cite{gomez_isogeometric_2010,nielsen_discretizations_2011,Bazilevs:2010:LES:1749635.1750210},
fluid-structure interaction problems 
\cite{bazilevs_isogeometric_2008,bazilevs_patient-specific_2009}. In addition, 
due to the ease of constructing high order continuous basis functions, IGA has been 
used with great success in solving PDEs that incorporate fourth order (or
higher)
derivatives of the field variable such as the Hill-Cahnard equation 
\cite{gomez_isogeometric_2008}, explicit gradient damage models \cite{verhoosel_isogeometric_2011-1} and gradient
elasticity \cite{fischer_isogeometric_2010}. 
The high order NURBS basis have also found potential application in the Kohn-Sham equation for electronic 
structure modeling of semiconducting materials \cite{Masud2012112}. We refer 
to \cite{nguyen_iga_review} for an overview of IGA and its implementation aspects.

In the context of fracture mechanics, IGA has been applied to fracture using the partition
of unity method (PUM) to capture two dimensional strong discontinuities and
crack tip singularities efficiently \cite{de_luycker_xfem_2011,ghorashi_extended_2012}.
In \cite{Tambat20121} an explicit isogeometric enrichment technique is proposed for modeling
material interfaces and cracks exactly. Note that this method
is contratry to PUM-based enrichment methods which define the cracks implicitly. 
A phase field model for dynamic fracture has been presented in \cite{Borden201277} where
adaptive refinement with T-splines provides an effective method for simulating fracture in three dimensions.
There are, however, only a few works on cohesive fracture in an IGA framework 
\cite{verhoosel_isogeometric_2011}. The method hinges on the ability to specify the continuity of NURBS/T-splines
through a process known as knot insertion. Highly accurate stress fields in cracked specimens were obtained with
coarse meshes.

Delamination or interfacial cracking between composite layers is unarguably 
one of the predominant modes of failure in
laminated composite. This failure mode has therefore been widely investigated both 
experimentally and numerically.
Delamination analyses have been traditionally performed using standard low order Lagrange finite
elements, see \eg \cite{Allix199561,Schellekens19931239,Crisfield,Krueger200125} and references therein. 
The two most popular computational methods for the analysis of delamination are the 
Virtual Crack Closure Technique (VCCT) \cite{rybicki1977,kruger2002} and interface elements 
with a cohesive law (also known as decohesion elements) 
\cite{Allix199561,Schellekens19931239,Crisfield}. The latter is adopted in this
contribution for it can deal with initiation and propagation of delamination in a unified theory. 
The Element Free Galerkin, which is a meshfree method, with the smooth 
moving least square basis was also adopted for delamination analysis \cite{Guiamatsia20092640}.
In order to alleviate the computational expense of cohesive elements, formulations with enrichment of the
FE basis was proposed in \cite{Samimi20112202,Guiamatsia20092616}. The extended finite element method 
(XFEM) \cite{mos_finite_1999} have been adopted for delamination studies \eg 
\cite{NME:NME907,doi:10.1142/S0219876206001181,CurielSosa2012788} which makes the pre-processing simple for
the delaminations can be arbitrarily located with respect to the FE mesh. The interaction between the 
delamination plane and the mesh is resolved during the solving step by using enrichment functions.
More recently, in \cite{Nguyen2013} high order B-splines cohesive FEM with $C^0$ 
continuity across element boundary
were utilized to efficiently model delamination of two dimensional (2D) composite specimens. In the referred
paper, it was shown that by using high order B-spline (order of up to 4) basis functions, relatively coarse
meshes can be used and 2D delamination benchmark tests such as the MMB were solved within 10 seconds on a laptop.

In this manuscript, prompted by our previous encouraging results reported on \cite{Nguyen2013} and
the work in \cite{verhoosel_isogeometric_2011},
we present an isogeometric framework for two  and three dimensional (2D/3D) delamination
analysis of laminated composites. Both the geometry and the displacement field are approximated using NURBS, therefore curved geometries are exactly represented.
We use knot insertion algorithm of NURBS to duplicate control
points along the delamination path. Meshes of zero-thickness interface elements can be
straightforwardly generated. The proposed ideas are implemented in our open source Matlab IGA code, 
MIGFEM\footnote{available for download at \url{https://sourceforge.net/projects/cmcodes/}}, 
described in \cite{nguyen_iga_review}. Several examples are provided including the mixed mode
bending test, a L-shaped curved composite specimen test \cite{lshape,Wimmer20082332}, 3D double cantilever beam
and a 3D singly curved thick-walled laminate. Moreover, isogeometric shell elements are used for the first 
time, at least to the authors' knowledge, to model delamination. Our findings are (i)
the proposed IGA-based framework reduces significantly the time being spent on the pre-processing step
to prepare FE models for delamination analyses and (ii) from the analysis perspective, 
the ultra smooth high order NURBS basis functions are able to produce highly accurate stress fields which is
very important in fracture modeling. The consequence is that relatively coarse meshes (compared to meshes of
lower order elements) can be adopted and thus the computational expense is reduced.

The remainder of the paper is organized as follows. Section \ref{sec:nurbs} briefly
presents NURBS curves, surfaces and solids.
Section \ref{sec:generation} is devoted to a discussion on knot insertion and
automatic generation of cohesive interface elements followed by finite element formulations
for solids with cohesive cracks given in Section \ref{sec:fem}. Numerical examples are given in
Section \ref{sec:examples}. Finally, Section \ref{sec:conclusions} ends the paper with some concluding remarks.

\section{NURBS curves, surfaces and solids}\label{sec:nurbs}

In this section, NURBS are briefly reviewed. 
We refer to the standard textbook \cite{piegl_book} for details.
A knot vector is a sequence in ascending order
of parameter values, written $\Xi=\{\xi_1,\xi_2,\ldots,\xi_{n+p+1}\}$
where $\xi_i$ is the \textit{i}th knot, $n$ is the number of basis functions and $p$ is 
the order of the B-spline basis. Open knots are used in this manuscript.

Given a knot vector $\Xi$, the B-spline basis functions are
defined recursively starting with the zeroth order basis
function ($p=0$) given by 

\begin{equation}
  N_{i,0}(\xi) = \begin{cases}
    1 & \textrm{if $ \xi_i \le \xi < \xi_{i+1}$}\\
    0 & \textrm{otherwise}
  \end{cases}
  \label{eq:basis-p0}
\end{equation}

\noindent and for a polynomial order $p \ge 1$

\begin{equation}
  N_{i,p}(\xi) = \frac{\xi-\xi_i}{\xi_{i+p}-\xi_i} N_{i,p-1}(\xi)
               + \frac{\xi_{i+p+1}-\xi}{\xi_{i+p+1}-\xi_{i+1}}
	       N_{i+1,p-1}(\xi)
  \label{eq:basis-p}
\end{equation}

\noindent This is referred to as the Cox-de Boor recursion formula.

\begin{figure}[h!]
  \centering 
  \includegraphics[width=0.6\textwidth]{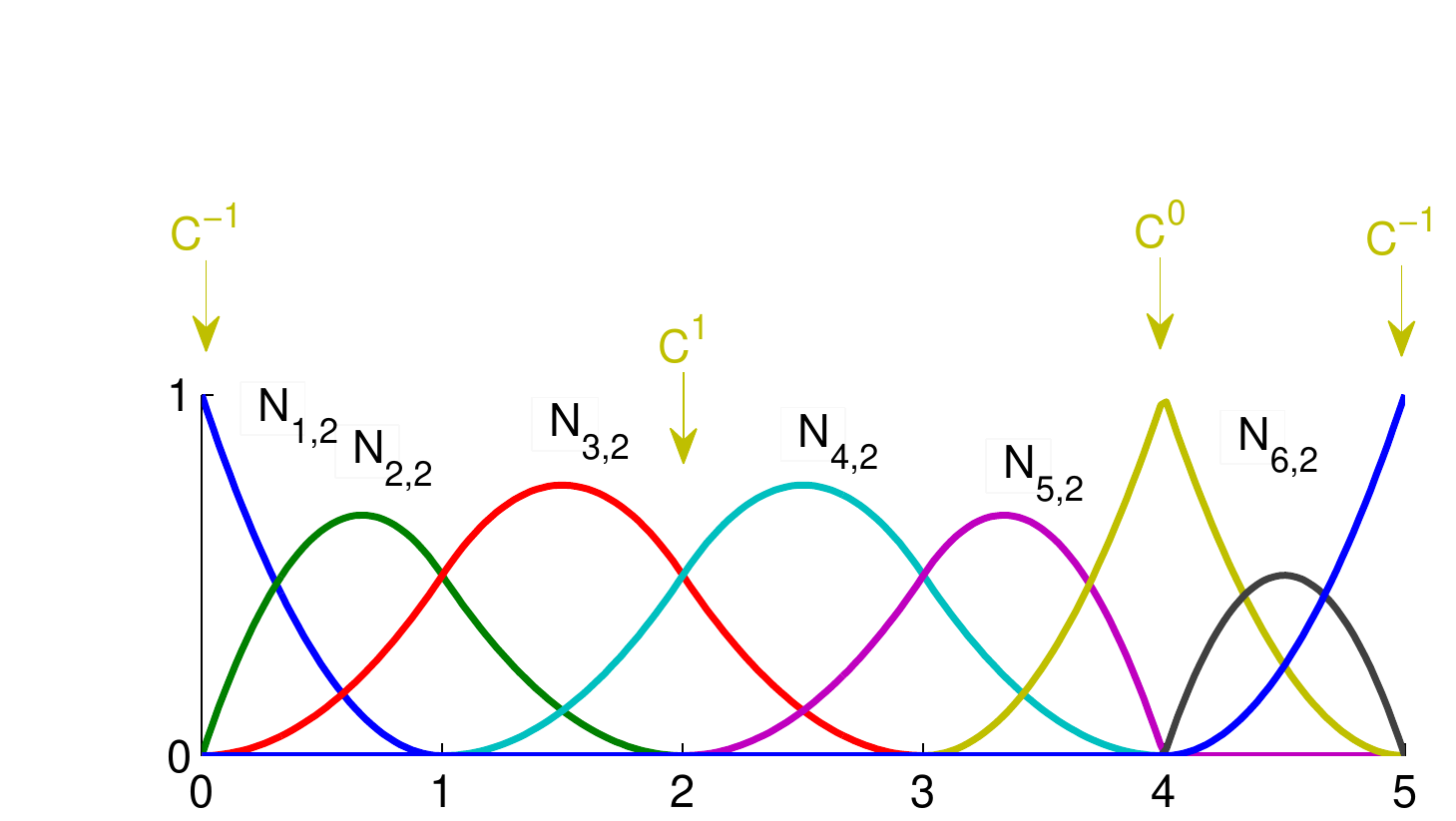}
  \caption{Quadratic ($p=2$) B-spline basis functions for an open non-uniform knot vector
  $\Xi=\{0,0,0,1,2,3,4,4,5,5,5\}$. Note the flexibility in the construction of
  basis functions with varying degrees of regularity.} 
  \label{fig:bspline-quad-open} 
\end{figure}

Figure \ref{fig:bspline-quad-open} illustrates some quadratic B-splines functions 
defined on an open non-uniform knot vector. 
Note that the basis functions are interpolatory at the
ends of the interval thanks to the use of open knot vectors 
and also at $\xi = 4$, the location of a repeated knot where
only $C^0$-continuity is attained. Elsewhere, the functions are $C^1$-continuous.
The ability to control continuity 
by means of knot insertion is particularly useful for modeling discontinuities
such as cracks or material interfaces as will be presented in this paper. In general,
in order to have a $C^{-1}$ continuity at a knot, its multiplicity must be $p+1$.

NURBS basis functions are defined as
\begin{equation}
	R_{i,p}(\xi) = \frac{N_{i,p}(\xi)w_i}{W(\xi)} =
	\frac{N_{i,p}(\xi)w_i}{\sum_{j=1}^{n}N_{j,p}(\xi)w_j},
  \label{eq:rational-basis}
\end{equation}
where $N_{i,p}(\xi)$ denotes the $i$th B-spline basis function of
order $p$ and $w_i$ are a set of $n$ positive weights. 
Selecting appropriate values for the $w_i$ permits the description of many
different types of curves including polynomials and circular arcs.
For the special case in which $w_i=c, i=1,2,\ldots,n$ the
NURBS basis reduces to the B-spline basis. Note that for simple geometries,
the weights can be defined analytically see \eg \cite{piegl_book}. For complex
geometries, they are obtained from CAD packages such as Rhino \cite{rhino}.

Given two knot vectors (one for each direction) $\Xi=\{\xi_1,\xi_2,\ldots,\xi_{n+p+1}\}$
and $\mathscr{H}=\{\eta_1,\eta_2,\ldots,\eta_{m+q+1}\}$ and a control net $\vm{B}_{i,j}\in \mathds{R}^d$,
a tensor-product NURBS surface is defined as 

\begin{equation}
	\vm{S}(\xi,\eta) =
	\sum_{i=1}^{n}\sum_{j=1}^{m}R_{i,j}^{p,q}(\xi,\eta)\vm{B}_{i,j}
  \label{eq:NURBS-surface1}
\end{equation}

\noindent where $R_{i,j}^{p,q}$ are given by

\begin{equation}
	R_{i,j}^{p,q}(\xi,\eta) = \frac{N_{i}(\xi) M_j(\eta) w_{i,j}}{
  \sum_{\hat{i}=1}^{n} \sum_{\hat{j}=1}^{m} N_{\hat{i}}(\xi) M_{\hat{j}}(\eta)
  w_{\hat{i},\hat{j}}}
  \label{eq:rational-basis2}
\end{equation}

\noindent In the same manner, NURBS solids are defined as 

\begin{equation}
	\vm{S}(\xi,\eta,\zeta) = \sum_{i=1}^{n}\sum_{j=1}^{m}\sum_{k=1}^l
	R_{i,j,k}^{p,q,r}(\xi,\eta,\zeta)\vm{B}_{i,k,j}
  \label{eq:NURBS-solid1}
\end{equation}

\noindent where $R_{i,j,k}^{p,q,r}$ are given by

\begin{equation}
  R_{i,j,k}^{p,q,r}(\xi,\eta,\zeta) = \frac{N_{i}(\xi) M_j(\eta) P_k(\zeta) w_{i,j,k}}{
  \sum_{\hat{i}=1}^{n} \sum_{\hat{j}=1}^{m}  \sum_{\hat{k}=1}^{l}N_{\hat{i}}(\xi) M_{\hat{j}}(\eta)
  P_{\hat{k}}(\zeta) w_{\hat{i},\hat{j},\hat{k}}}
  \label{eq:rational-basis3}
\end{equation}
\noindent Derivatives of the B-splines and NURBS basis functions can be find elsewhere \eg
\cite{hughes_isogeometric_2005,cottrel_book_2009}.

\section{Automatic generation of cohesive elements}\label{sec:generation}

\subsection{Knot insertion}\label{sec:knotInsertion}

It should be emphasized that knot insertion does not change the B-spline
curves or surfaces geometrically but a direct influence on the continuity of
the approximation where knots are repeated. Let us consider a knot vector defined by
$\Xi=\{\xi_1,\xi_2,\ldots,\xi_{n+p+1}\}$ with  the corresponding control net 
denoted by $\vm{B}$. A new extended knot vector given by
$\bar{\Xi}=\{\bar{\xi}_1=\xi_1,\bar{\xi}_2,\ldots,\bar{\xi}_{n+m+p+1}=
\xi_{n+p+1}\}$ is formed where $m$ knots are added. The $n+m$ new control points
$\bar{\vm{B}}_i$ are formed from the original control points by

\begin{equation}
  \bar{\vm{B}}_i = \alpha_i \vm{B}_i + (1-\alpha_i)\vm{B}_{i-1}
  \label{eq:new-control-pnt}
\end{equation}
\noindent where 
\begin{equation}
  \alpha_i = \begin{cases}
    1 & 1 \le i \le k-p,\\
    \D\frac{\bar{\xi}-\xi_i}{\xi_{i+p}-\xi_i} & k-p +1 \le i \le k\\
    0 & k+1 \le i \le n+p+2
  \end{cases}
  \label{eq:alpha-i}
\end{equation}

\begin{figure}[htbp]
  \centering 
  \includegraphics[width=0.33\textwidth]{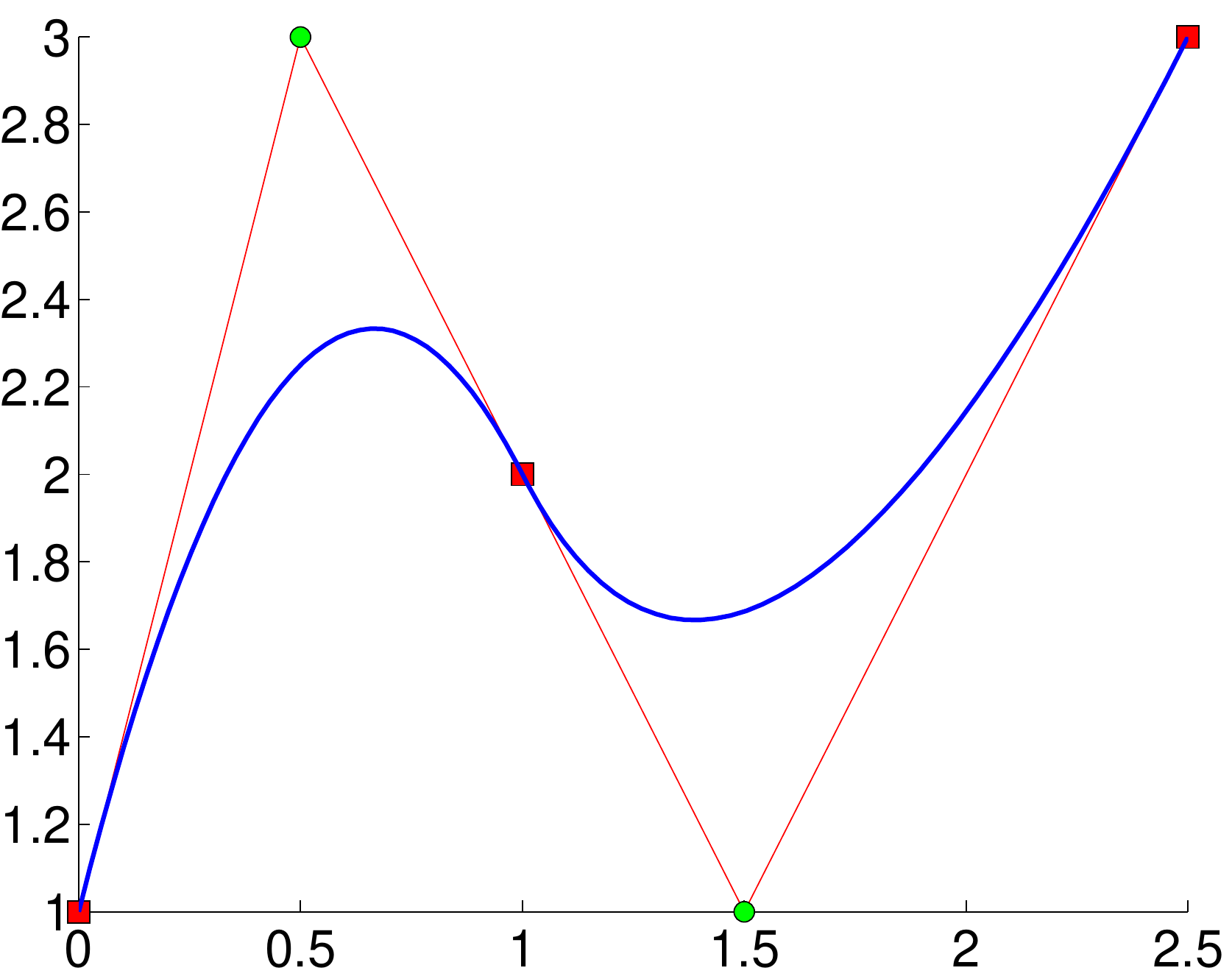}
  \includegraphics[width=0.33\textwidth]{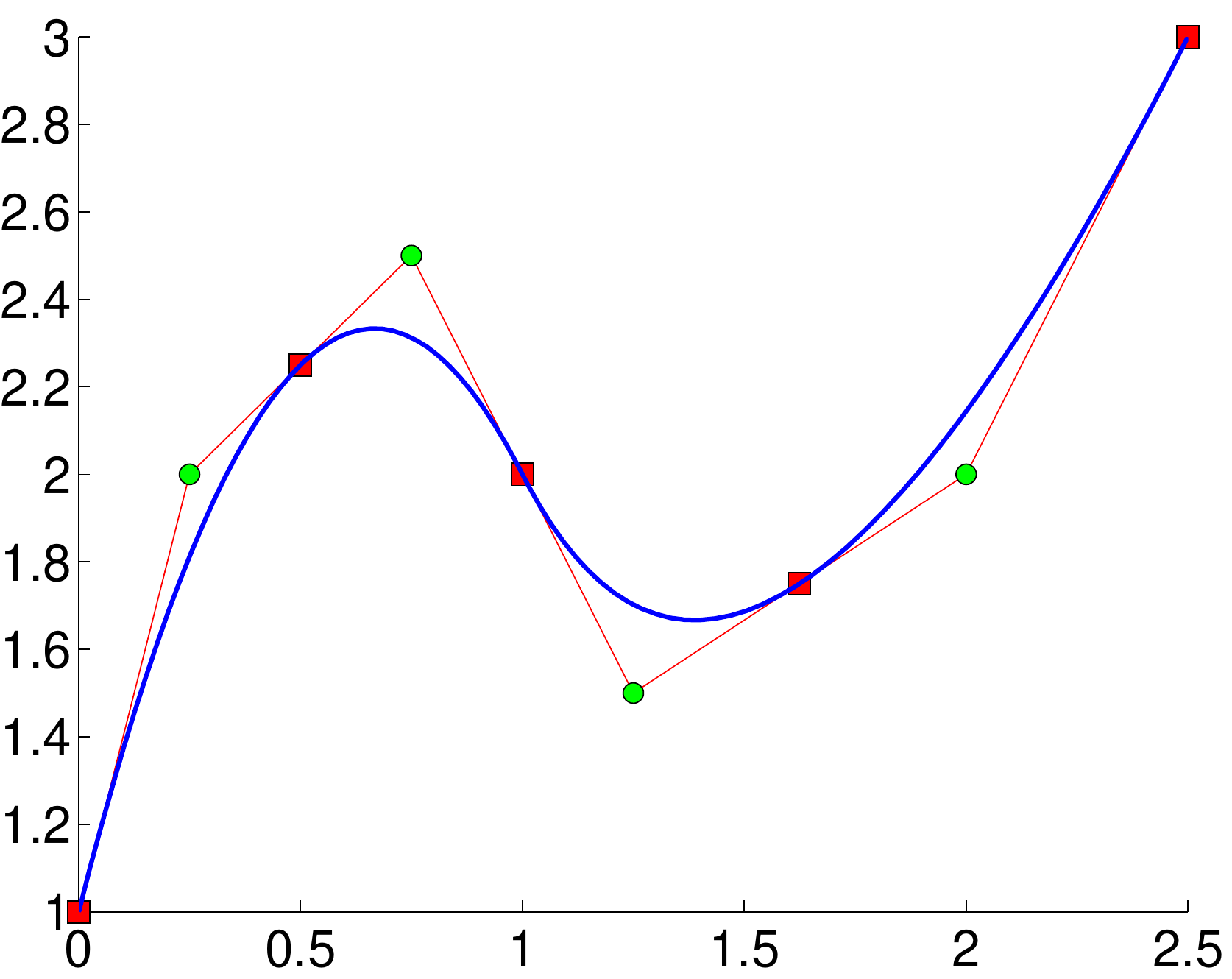}
  \caption{Knot insertion on a quadratic B-spline curve. The curve is not
  changed geometrically. Control points are denoted by filled green circles.
  Points corresponding to the knot values are denoted by red circles. These
  points divide the curve into segments or elements from an analysis standpoint.}
  \label{fig:refine1d} 
\end{figure}

Considering a quadratic B-spline curve with 
knot vector $\Xi=\{0,0,0,0.5, 1, 1, 1\}$ and control points as shown in Fig.
\ref{fig:refine1d} (left). On the right of the same figure, two new knots $\xi=0.25$ and
$\xi=0.75$ were added. Consequently, two new control points were formed.
Although the curve is not changed geometrically and parametrically, the basis functions are now richer and may be more suitable for the purpose of analysis.

\begin{figure}[h!]
  \centering 
  \subfloat[$\Xi=\{0,0,0,1,1,1\}$]{\includegraphics[width=0.3\textwidth]{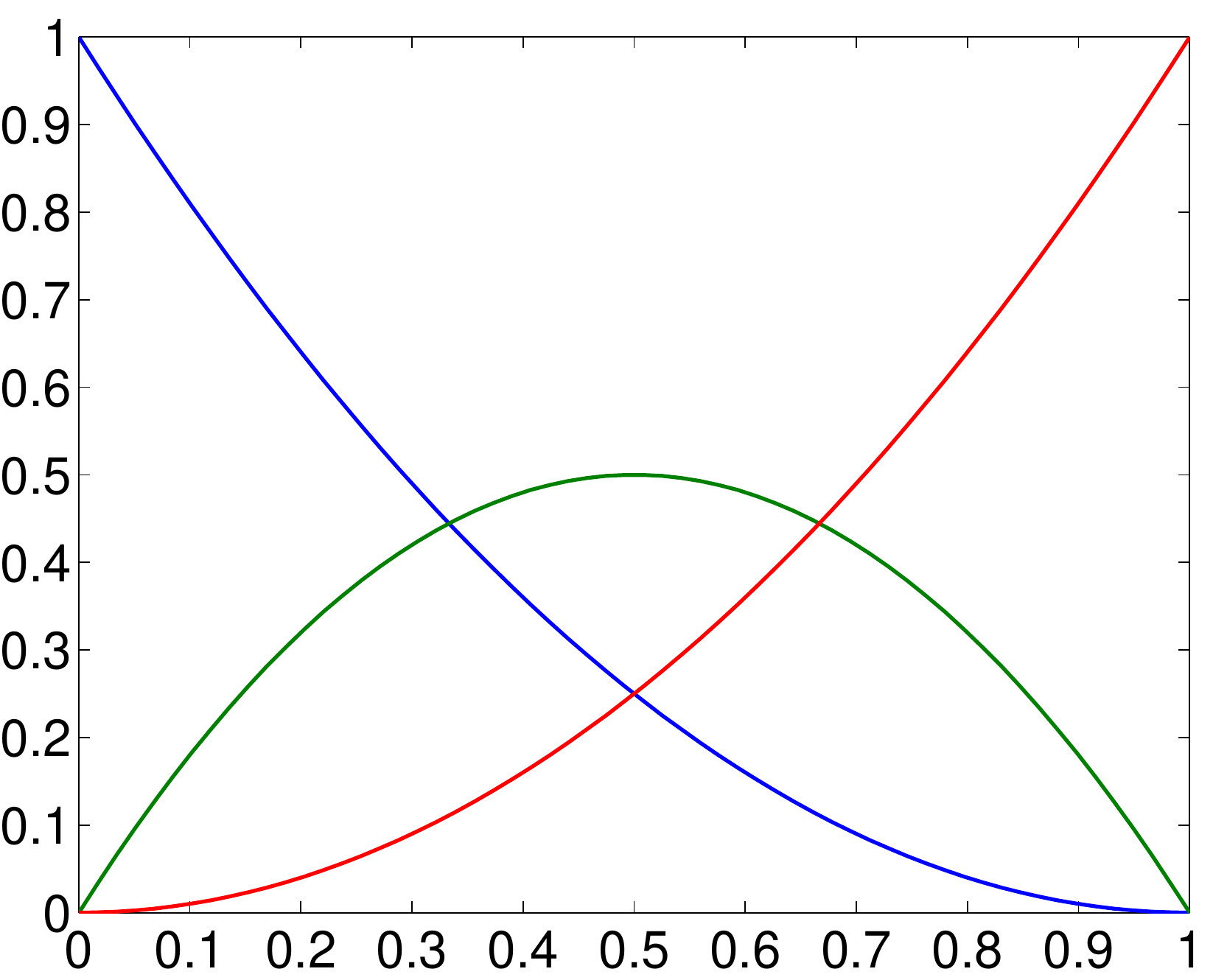}}\;
  \subfloat[$\Xi'=\{0,0,0,0.5,0.5,0.5,1,1,1\}$]{\includegraphics[width=0.3\textwidth]{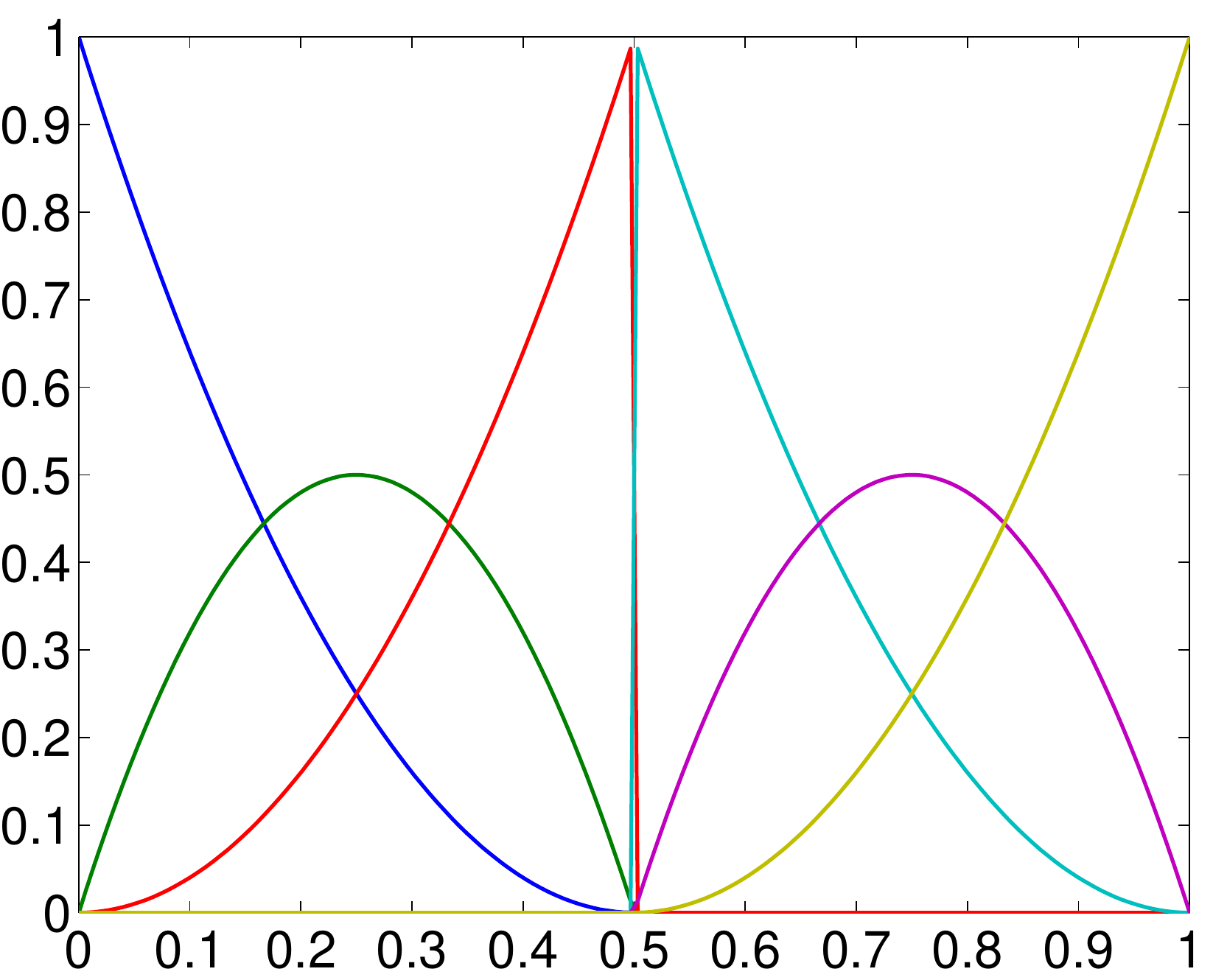}}\\
  \subfloat[B-spline curve after knot insertion]{\includegraphics[width=0.3\textwidth]{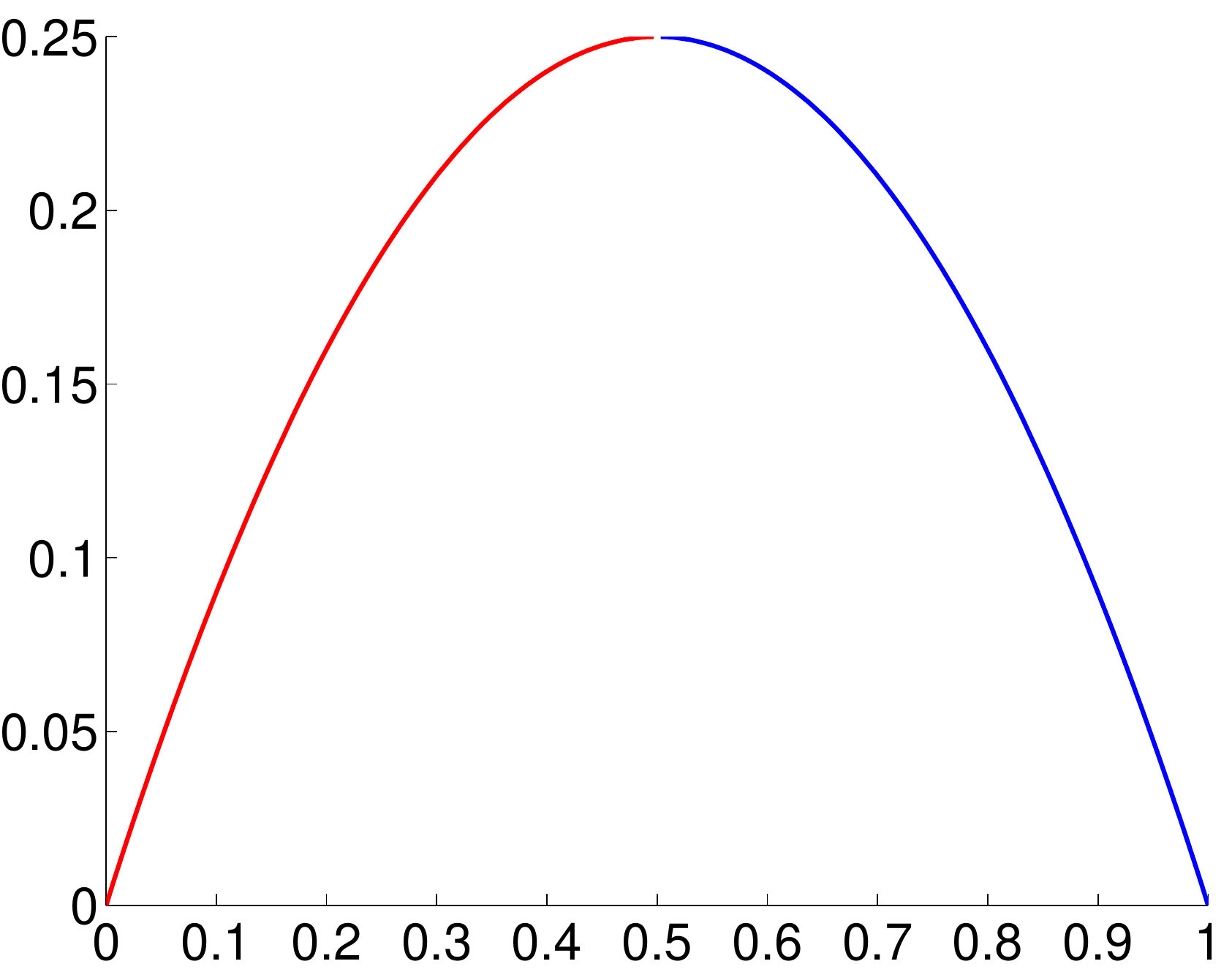}}
  \subfloat[B-spline curve with $\vm{B}'_4$ moved slightly the original
  position]{\includegraphics[width=0.3\textwidth]{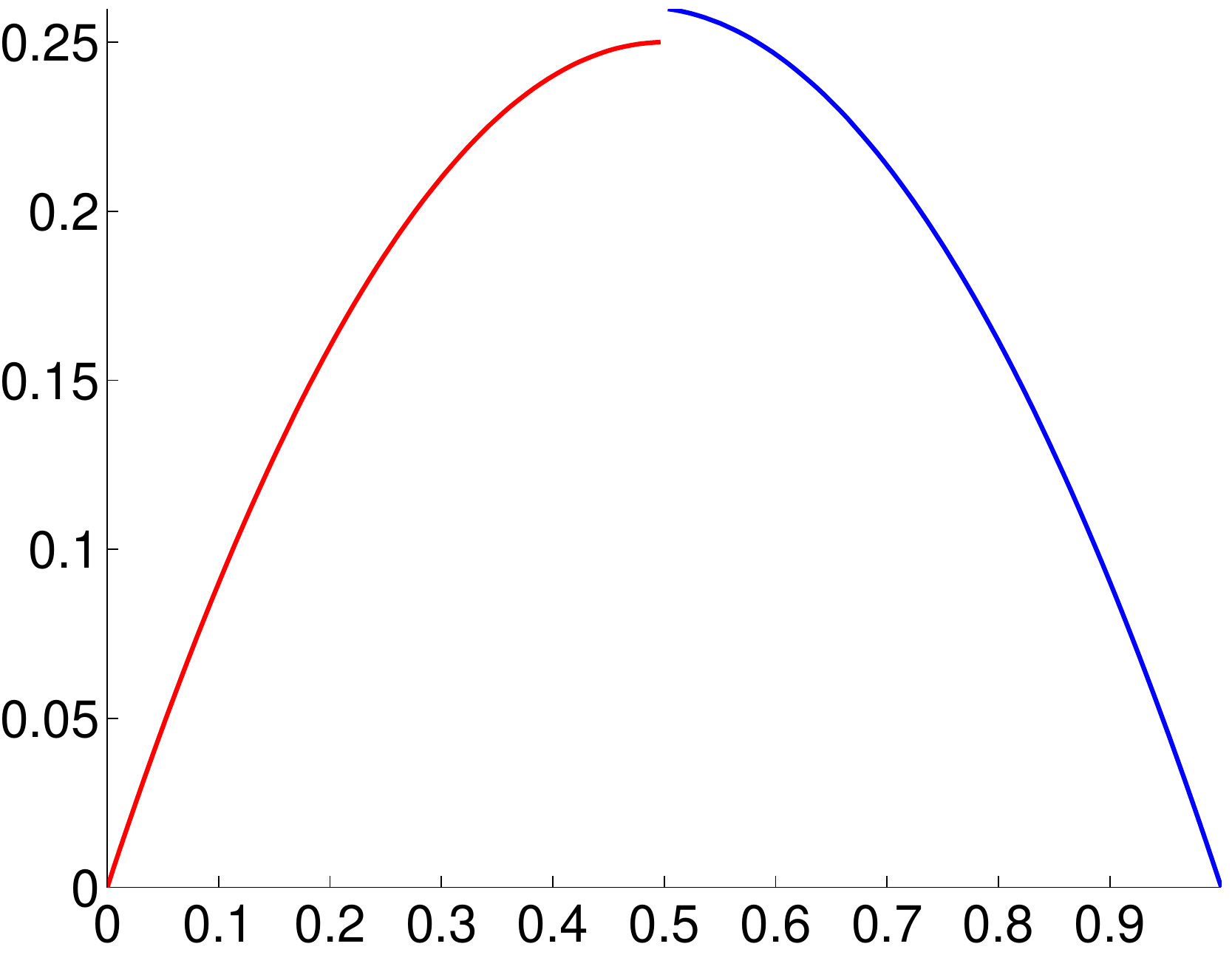}}
  \caption{$p+1$ times knot insertion for a quadratic B-spline curve to
  introduce a $C^{-1}$ discontinuity at $\xi=0.5$.}
  \label{fig:discontinuity1d} 
\end{figure}

Let us now consider a quadratic B-spline defined using $\Xi=[0,0,0,1,1,1]$. The
three basis functions for this curve are given in Fig. (\ref{fig:discontinuity1d}a). Now
suppose that we need to have a discontinuity at $\xi=0.5$. This can be
achieved by inserting a new knot $\bar{\xi}=0.5$ three ($=p+1$) times. 
The new knot vector is then given by $\Xi'=[0,0,0,0.5,0.5,0.5,1,1,1]$ and the new
basis functions are shown  in Fig. (\ref{fig:discontinuity1d}b). 
Let us build a B-spline curve with the control net defined by $\vm{B}$ as shown in
Eq. (\ref{eq:discontinuity1}). The new control net that is defined by
$\vm{B}'$ is also given in Eq. (\ref{eq:discontinuity1}).

\begin{equation}
	\vm{B} = \begin{bmatrix}
		0.0 & 0.0 \\
		0.5 & 0.5 \\
		1.0 & 0.0
	\end{bmatrix}, \quad \vm{B}'=\begin{bmatrix}
		0.00 & 0.00 \\
		0.25 & 0.25 \\
		0.50 & 0.25 \\
		0.50 & 0.25 \\
		0.75 & 0.25 \\
		1.00 & 0.00 \\
	\end{bmatrix}.
	\label{eq:discontinuity1}
\end{equation}

\noindent where it should be noted that $\vm{B}'_3=\vm{B}'_4$. The B-spline curve corresponds to
the original and new basis is the same and given in Fig.
(\ref{fig:discontinuity1d}c). Imagine now that point $\vm{B}'_4$ slightly
moves vertically, the resulting B-spline curve with a strong discontinuity at
$x=0.5$ is plotted in Fig. (\ref{fig:discontinuity1d}d). This technique of
inserting knot values $p+1$ times was used in \cite{verhoosel_isogeometric_2011} to
model the decohesion of material interfaces. The application of this method in
two/three dimensions resemble the usage of zero-thickness interface elements by
doubling nodes in the FE framework.

\begin{figure}[h!]
  \centering 
  \psfrag{a}{(a) $\Xi=\mathscr{H}=\{0,0,1,1\}$}
  \psfrag{b}{(b) $\overline{\mathscr{H}}=\{0.5,0.5,0.5\}$}
  \psfrag{c}{(c) Duplicated nodes are shown}
  \psfrag{d}{(d) Refined mesh}
  \includegraphics[width=0.95\textwidth]{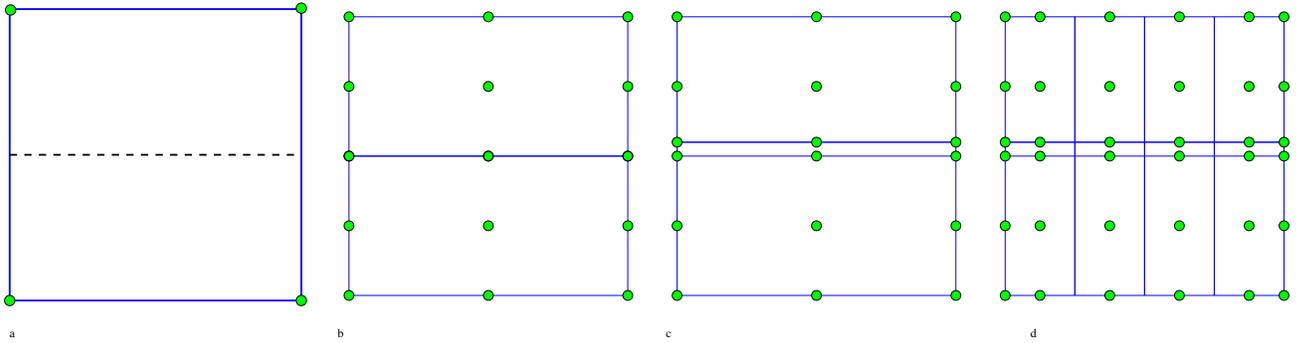}
  \caption{Example of introducing a horizontal discontinuity in a NURBS surface.}
  \label{fig:square}
\end{figure}

\begin{figure}[h!]
  \centering 
  \psfrag{n}{s}
  \includegraphics[width=0.5\textwidth]{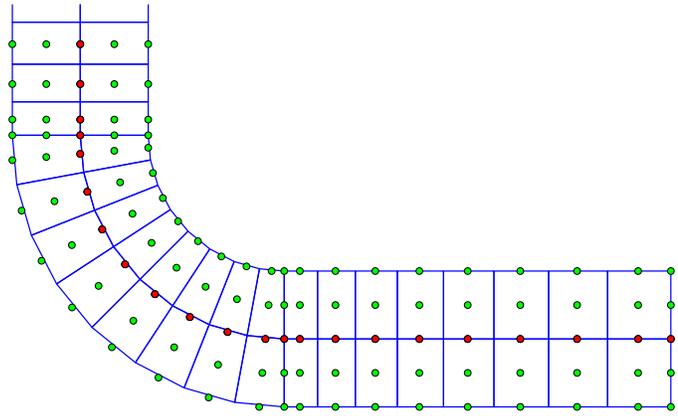}
  \caption{L-shaped composite sample of two plies with a fillet modeled with a bi-quadratic NURBS: 
         red circles denote duplicated nodes.
         For this case, it suffices to find the index of node S--the first node on the discontinuity curve.
         By virtue of the tensor-product nature of NURBS, the indices of other discontinuity nodes can then be
         found with ease.}
  \label{fig:lshape}
\end{figure}

We demonstrate the technique to generate a discontinuity into a NURBS surface by a simple example.
The studied surface is a square of $10\times10$ and suppose that one needs a horizontal discontinuity
line in the middle of the square as shown in Fig. (\ref{fig:square}a). The coarsest mesh consists of 
one single bi-linear NURBS element with $\Xi=\mathscr{H}=\{0,0,1,1\}$ and $p=q=1$. To insert the desired
discontinuity, the following steps are performed:
(1) perform order elevation to $p=q=2$;
(2) perform knot insertion for $\mathscr{H}$, the new knot is  $\mathscr{H}=\{0,0,0,0.5,0.5,0.5,1,1,1\}$ 
(Fig. (\ref{fig:square}b)); and
(3) perform knot insertion to refine the mesh if needed.
In Fig. (\ref{fig:square}c,d) the duplicated control points were moved upward to show the effect of discontinuity.
In order to use these duplicated nodes in a FE context, one can put springs connecting each pair of nodes or
employ zero-thickness interface elements. In this manuscript the latter is used. With a small amount of effort, the connectivity matrix for the interface elements can be constructed using a simple Matlab code as 
given in Listing \ref{list-1D}. It is obvious that, due to the simplification made in line 2 of 
Listing \ref{list-1D}, this code snippet applies only for a horizontal/vertical discontinuity line. However,
it is straightforward to extend this template code to general cases by changing line 2. Such refinements
are certainly problem dependent and hence not provided here. We refer to 
Fig. (\ref{fig:lshape}) for one example of a curved composite panel made of two plies.

\begin{snippet1}[caption={Matlab code to build the element connectivity for 1D interface elements}, 
   label={list-1D},framerule=1pt]
    [ielements] = buildIGA1DMesh (uKnot,p);
    delaminationNodes  =  find(abs(controlPts(:,2) - 5  ) <1e-10);
    mm                 = 0.5*length(delaminationNodes);
    lowerNodes         = delaminationNodes(1:mm);
    upperNodes         = delaminationNodes(mm+1:end);
    
    iElements   = zeros(noElemsU,2*(p+1));
    
    for i=1:noElemsU
        sctr = ielements(i,:);
        iElements(i,1:p+1)   = lowerNodes(sctr);
        iElements(i,p+2:end) = upperNodes(sctr);
    end                
\end{snippet1}

\begin{snippet1}[caption={Matlab code to build the element connectivity for 2D interface elements}, 
   label={list-2d},framerule=1pt]
    delaminationNodes  =  find(abs(controlPts(:,3) -b/2  ) <1e-10);
    mm                 = 0.5*length(delaminationNodes);
    lowerNodes         = delaminationNodes(1:mm);
    upperNodes         = delaminationNodes(mm+1:end);
    % noElemsU=number of elements along X-dir
    iElements   = zeros(noElemsU*noElemsV,2*(p+1)*(q+1));
    iElementS   = generateIGA2DMesh (uKnot,vKnot,noPtsX,noPtsY,p,q);
    for e=1:noElemsU*noElemsV
        iElements(e,1:(p+1)*(q+1))     = lowerNodes(iElementS(e,:));
        iElements(e,(p+1)*(q+1)+1:end) = upperNodes(iElementS(e,:));
    end   
\end{snippet1}

\begin{figure}[h!]
  \centering 
  \psfrag{x}{$X$} \psfrag{y}{$Y$} \psfrag{z}{$Z$}
  \psfrag{dis}{discontinuity surface}
  \includegraphics[width=0.2\textwidth]{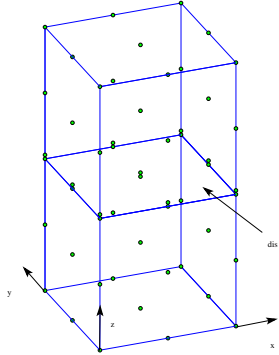}
  \caption{A 3D bar with a discontinuity surface in the middle: modeled by a tri-quadratic NURBS solid.}
  \label{fig:cube}
\end{figure}

The technique introduced so far can be straightforwardly extended to three dimensions, see 
Listing \ref{list-2d} and Fig. (\ref{fig:cube}) for an example.
The discontinuity surface lies in the $X-Y$ plane. Line 7 of this Listing builds the
element connectivity array for a 2D NURBS mesh, we refer to \cite{nguyen_iga_review} for a detailed
description of these Matlab functions.
These pre-processing techniques are implemented in our open source Matlab IGA code named MIGFEM,
desribed in \cite{nguyen_iga_review},      
which is available at \url{https://sourceforge.net/projects/cmcodes/}.
In order to support IGA codes which are based on the B{\'e}zier extraction 
\cite{borden_isogeometric_2011,scott_isogeometric_2011}, see also Section \ref{sec:implementation}, 
MIGFEM computes the 1D, 2D and 3D B{\'e}zier extractors. In summary the pre-processing code 
writes to a file with (1) coordinates of control points (including duplicated ones), 
(2) connectivity of continuum elements, (3) connectivity of interface elements, (4)
2D/3D B{\'e}zier extractors for continuum elements and (5) 1D/2D extractors for interface
elements. It should be emphasized that inserting interface elements into a Lagrange FE mesh 
is a time-consuming task even with commercial FE packages. Due to that fact, a free mesh generator 
for cohesive modeling was developed by the first author and presented in \cite{PhuInterface}.

\begin{rmk}
In the proposed framework, interface elements are inserted \textit{a priori}, therefore delaminations
only grow along predefined paths. For laminates built up by plies of unidirectional fiber reinforced
composites, the fracture toughness of the plies is much greater than the fracture toughness of the ply 
interfaces. Therefore, delaminations only grow along the ply interfaces which are known \textit{a priori}.
And that justifies our assumption.
\end{rmk}

\section{Finite element formulation}\label{sec:fem}

\subsection{Isogeometric analysis}

According to the IGA the field variable (which is, in this paper, the displacement field) is 
approximated by the same B-spline/NURBS basis functions used to exactly represent the geometry. 
Therefore, in an IGA context, one writes for the geometry and displacement field, respectively

\begin{subequations}
\begin{align}
\vm{x} &= N_I (\bsym{\xi})\vm{x}_I\\
u_{i} &= N_I(\bsym{\xi}) u_{iI}
\end{align}
\label{eq:isoparametric}
\end{subequations}

\noindent where $\vm{x}_I$ are the nodal coordinates, $u_{iI}$ is the $i$ ($i=1,2,3$) 
component of the displacement at node/control point $I$ and $N_I$ denotes the shape functions which are
the B-spline/NURBS basis functions described in Section \ref{sec:nurbs}.

\begin{figure}[htbp]
  \centering 
  \psfrag{x}{x}
  \psfrag{y}{y}
  \psfrag{parent}{$\square$}
  \psfrag{xib}{$\bar{\xi}$}
  \psfrag{xi}{$\xi$}
  \psfrag{eta}{$\eta$}
  \psfrag{etab}{$\bar{\eta}$}
  \psfrag{xi1}{$\xi_i$}
  \psfrag{xi2}{$\xi_{i+1}$}
  \psfrag{eta1}{$\eta_{j}$}
  \psfrag{eta2}{$\eta_{j+1}$}
  \psfrag{c1}{$(-1,-1)$}
  \psfrag{c2}{$(1,1)$}
  \psfrag{para}{$\hat{\Omega}_e$}
  \psfrag{phy}{$\Omega_e$}
  \psfrag{xi-line}{$\xi$-lines}
  \psfrag{eta-line}{$\eta$-lines}
  \psfrag{image}{images of $\xi$ and $\eta$-lines}
  \psfrag{a}{Physical domain}\psfrag{b}{Parametric domain}\psfrag{c}{Parent
  domain}
  \includegraphics[width=0.7\textwidth]{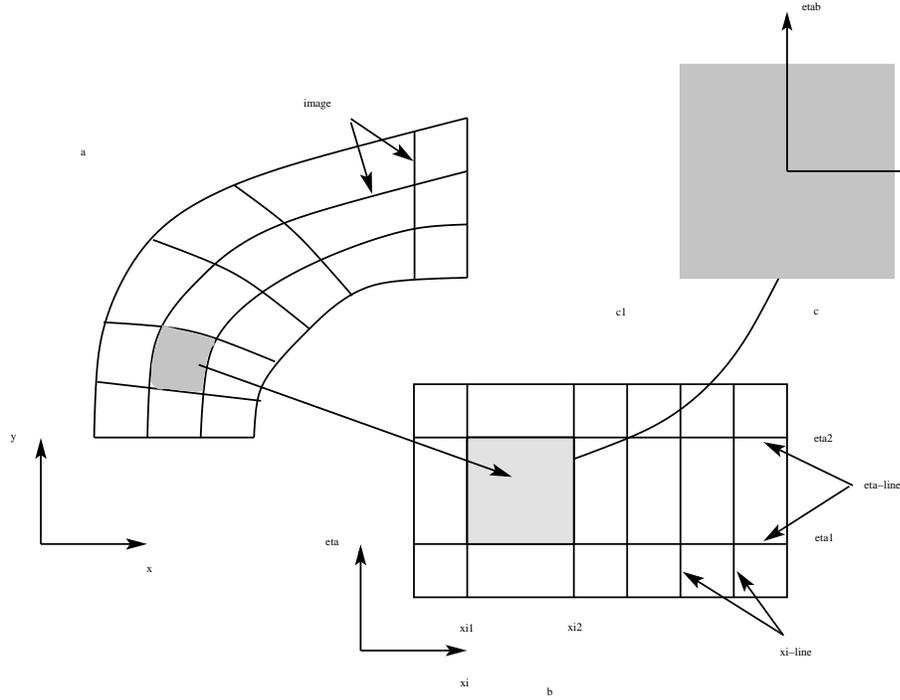}
  \caption{Definition of domains used for integration in isogeometric
    analysis. Elements are defined in the parametric space as non-zero knot spans, $[\xi_i,\xi_{i+1}]\times
       [\eta_j,\eta_{j+1}]$ and elements in the physical space are images of their parametric counterparts.} 
  \label{fig:integration} 
\end{figure}

Elements are defined as non-zero knot spans, see Fig. (\ref{fig:integration}), which are elements
in the parameter space (denoted by $\hat{\Omega}_e$). Their images in the physical space 
obtained via the mapping, see Eq. (\ref{eq:isoparametric}), are called elements in the 
physical space (denoted by $\Omega_e$) that resemble the familiar
Lagrange elements. From our experiences, it is beneficial to work with elements in the parameter space.
Numerical integration is also performed on a parent domain as in Lagrange FEs.

\subsection{FE discrete equations}

The semi-discrete equation for a solid with cohesive cracks is given by

\begin{equation}\label{semi-discrete}
 \vm{M} \vm{a} = \vm{f}^{\textrm{ext}}  - \vm{f}^{\textrm{int}} - \vm{f}^{\textrm{coh}}
\end{equation}

\noindent where $\vm{a}$ is the acceleration vector, $ \vm{M} $ denotes the consistent mass matrix, 
$ \vm{f}^{\textrm{ext}}$ is the external force vector
, the internal force vector is denoted as $\vm{f}^{\textrm{int}}$ and 
the cohesive force vector $\vm{f}^{\textrm{coh}}$. The elemental mass matrix, external and internal force vectors
are computed from contributions of continuum elements and given by

\begin{align}\label{my}
 \vm{M}_{e}  &=  \int_{\Omega_e} \rho \vm{N}^{\mathrm{T}} \vm{N} \mathrm{d} \Omega_{e}\\
 \vm{f}^{\textrm{int}}_{e}  &= \int_{\Omega_e} \vm{B}^{\mathrm{T}} \bsym{\sigma} \mathrm{d} \Omega_{e}\\
 \vm{f}^{\textrm{ext}}_{e} &= \int_{\Omega_e} \rho \vm{N}^{\mathrm{T}} \vm{b} \mathrm{d} \Omega_{e} + \int_{\Gamma_t^{e}} \vm{N}^{\mathrm{T}} \bar{\vm{t}}  \mathrm{d} \Gamma_t^{e}
\end{align}

\noindent where $\rho$ is the density, $\Omega_e$ is the element domain, $\Gamma_t^{e}$ is the
element boundary that overlaps with the Neumann boundary, $\vm{b}$ and 
$\bar{\vm{t}}$ are the body forces and traction vector, respectively. The shape function matrix
and the strain-displacement matrix are denoted by $\vm{N}$ and $\vm{B}$; $\bsym{\sigma}$ is the
Cauchy stress vector.

The cohesive force vector is computed by assembling the contribution
of all interface elements. It is given by for an interface element $ie$

\begin{equation}\label{fcoh}
\begin{split}
\vm{f}^{\textrm{coh}}_{ie,+} &= 
\int_{\Gamma} \vm{N}^{\mathrm{T}}_\text{int} \vm{t}^{\mathrm{c}}  \mathrm{d} \Gamma \\
\vm{f}^{\textrm{coh}}_{ie,-} &= -
\int_{\Gamma} \vm{N}^{\mathrm{T}}_\text{int} \vm{t}^{\mathrm{c}}  \mathrm{d} \Gamma 
\end{split}
\end{equation}

\noindent in which $\vm{t}^c$ denotes the cohesive traction, $\vm{N}_\text{int}$ represents
the shape function matrix of interface elements. The subscripts +/- denote the upper and lower
faces of the interface element.

\begin{figure}[h!]
  \centering 
  \psfrag{solid}{solid elements}
  \psfrag{interface}{$C^{p-1}$ interface elements}
  \psfrag{interface1}{$C^0$ interface elements}
  \psfrag{n1}{$N_1(\xi)$} \psfrag{n2}{$N_2(\xi)$} \psfrag{n3}{$N_3(\xi)$} \psfrag{n4}{$N_4(\xi)$}
  \psfrag{1}{1} \psfrag{2}{2} \psfrag{3}{3} \psfrag{4}{4}
  \psfrag{5}{5} \psfrag{6}{6} \psfrag{7}{7} \psfrag{8}{8}
  \psfrag{9}{9} \psfrag{10}{10} 
  \includegraphics[width=0.7\textwidth]{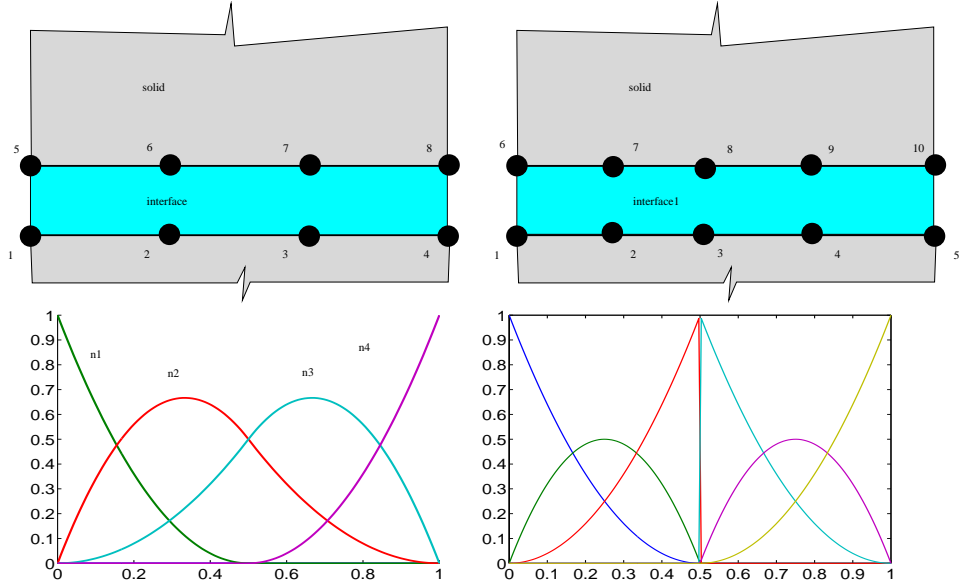}
  \caption{Illustration of $C^{p-1}$ NURBS interface elements (left) and $C^0$ NURBS interface elements (right):
  For $C^{p-1}$ elements, the connectivity of the first element is $[1,2,3,5,6,7]$
  while the connectivity of the second element is $[2,3,4,6,7,8]$. }
  \label{fig:interface}
\end{figure}

The displacement of the upper and lower faces of an interface element, let say the first
element in Fig. (\ref{fig:interface})-left read

\begin{equation}\label{lower-upper}
\begin{split}
\vm{u}^+ &= N_1(\xi) \vm{u}_5 + N_2(\xi) \vm{u}_6 + N_3(\xi) \vm{u}_7 \\
\vm{u}^- &= N_1(\xi) \vm{u}_1 + N_2(\xi) \vm{u}_2 + N_3(\xi) \vm{u}_3 \\
\end{split}
\end{equation}

\noindent with $N_I$ ($I=1,2,3,4$) are the quadratic NURBS shape functions.
Figure \ref{fig:interface} also explains the difference between $C^{p-1}$ and $C^0$ high order
elements--for the same number of elements, $C^{p-1}$ meshes have less nodes. We refer to
\cite{hughes_isogeometric_2005} for more information on this issue.
The latter was used in \cite{Nguyen2013} with B-spline basis for 2D delamination analysis.

Having defined the displacement of the upper and lower faces of the
interface, it is able to compute the displacement jump as

\begin{equation}\label{jump}
 \llbracket \vm{u}(\vm{x}) \rrbracket \equiv \vm{u}^+ - \vm{u}^-= \vm{N}_\text{int} ( \vm{u}^+ - \vm{u}^- )
\end{equation}
\noindent where 
\begin{equation}
\vm{N}_\text{int}=\begin{bmatrix}
N_1 & 0 & N_2 & 0 & N_3 & 0\\ 0 & N_1 & 0 & N_2 & 0 & N_3
\end{bmatrix},\quad
\vm{u}^+=\begin{bmatrix}
u_{x5} & u_{y5} & u_{x6} & u_{y6} & u_{x7} & u_{y7}
\end{bmatrix}\trans
\end{equation}

\noindent The displacement jump will be inserted into a cohesive law (or traction-separation law)
to compute the corresponding traction $\vm{t}^c$. We refer to 
\cite{Allix199561,Schellekens19931239,Crisfield} and references therein for other aspects of 
interface cohesive elements. The implementation for three dimensional problems \ie 2D interface
elements is straightforward, for example in Eq. (\ref{lower-upper}), instead of using univariate 
NURBS basis one uses bivariate basis $N_I(\xi,\eta)$.

\subsection{Cohesive laws}

In this work, we adopt the damage-based bilinear cohesive law developed in \cite{Camanho01082003,Turon20061072}.
This is a cohesive law in which the fracture toughness is a phenomenological function, rather than a material constant, of mode mixity as formulated by Benzeggagh and Kenane \cite{Benzeggagh1996439}.
Herein we briefly recall the cohesive law of which implementation details can be found in
\cite{vanderMeer2010719}. Denoting $d$ as the damage variable ($0\le d \le 1$), the cohesive
law reads in the local coordinate system attached to the interface elements

\begin{equation}
\vm{t}^c_l= (1-d)K \llbracket \vm{u} \rrbracket_l
\end{equation}
\noindent where $K$ is the dummy stiffness. The damage variable $d$ is a function of the
equivalent displacement jump, the onset $\jump{u}_\text{eq}^0$
and the propagation equivalent displacement jump  $\jump{u}_\text{eq}^f$.
The onset $\jump{u}_\text{eq}^0$ is a function of $K$, the mode mixity and the normal and shear strength 
$\tau^0_1$ and $\tau^0_3$. The propagation displacement jump $\jump{u}_\text{eq}^f$ is a function of
$\jump{u}_\text{eq}^0$, mode I and II fracture toughness $G_{Ic}$, $G_{IIc}$, the mode mixity and $\eta$
which is a curve fitting value for fracture toughness tests performed by  
Benzeggagh and Kenane \cite{Benzeggagh1996439}.

\subsection{Numerical integration}

In this manuscript, full Gaussian integration schemes are used. Precisely,
for 2D solid elements of order $p\times q$, a $(p+1)\times(q+1)$ Gauss quadrature rule is 
adopted and for cohesive elements of order $p$, a $(p+1)$ Gauss scheme is utilized. A similar rule
was used for 3D solid elements and 2D cohesive elements.

\subsection{Implementation aspects}\label{sec:implementation}

There are at least two approaches to incorporating IGA into existing FE codes--with and without
using the B{\'e}zier extraction. The former, which relies on the B{\'e}zier decomposition technique, 
was developed in 
\cite{borden_isogeometric_2011,scott_isogeometric_2011} and provides data structures (the so-called
B{\'e}zier extractor sparse matrices) that facilitate the implementation of IGA in existing FE codes.
Precisely, the shape functions of IGA elements are the Bernstein polynomials (defined on the standard
parent element) multiplied by the extractors. 
We refer to \cite{nguyen_iga_review} for a discussion on both techniques. 

\begin{figure}[htbp]
         \centering
         \psfrag{et}[c]{$\eta$ lines}
         \psfrag{xi}[c]{$\xi$ lines}
         \psfrag{vm}[c]{visualization Q4 mesh}
         \psfrag{gp}{Gauss points}
         \psfrag{vn}{visualization nodes}
         \includegraphics[width=0.6\textwidth]{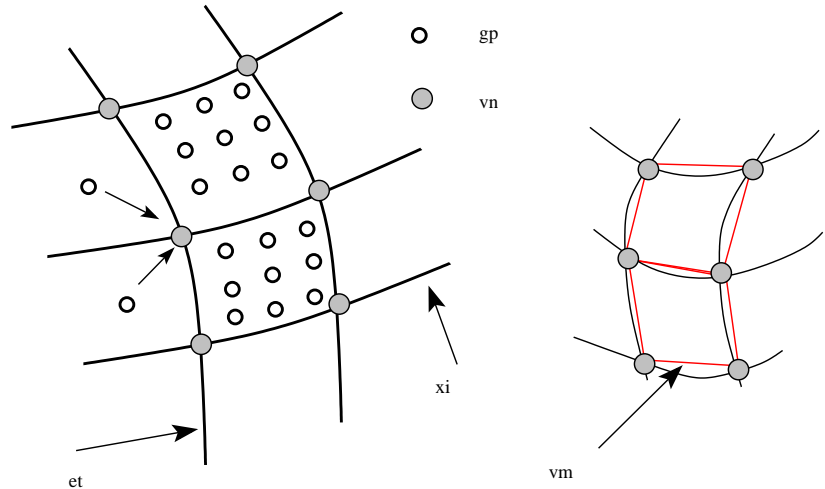}
        \caption{Post-processing in Isogeometric Analysis.}
        \label{fig:visualization}
\end{figure}

For curved geometries, the post-processing of IGA is more involved than Lagrange FEs due to two reasons
(1) some control points locate outside the physical domain (hence the computed displacements
at control points are not nodal values) and (2) existing post-processing techniques cannot be applied
directly to NURBS meshes. Interested reader can refer to \cite{nguyen_iga_review} for a discussion on
some post-processing techniques for IGA. For completeness we discuss briefly one technique here for 2D
problems. First, a visualization mesh which consists of four-noded quadrilateral elements is constructed.
The nodes of this mesh are the intersections of the $\xi$ and $\eta$ knot lines in the physical space.
We then extrapolate the quantities at Gauss points to these nodes and perform nodal averaging 
if necessary. Figure \ref{fig:visualization} summarizes the idea.

\section{Examples}\label{sec:examples}

Since we are introducing a computational framework for delamination analyses rather than a detailed
study of the delamination behaviour of composite materials, intralaminar damage (matrix cracking and
fiber damage) is not taken into account leading to an orthotropic elastic behaviour assumption for the plies.
Note that matrix cracking can however be efficiently modeled using extended finite elements as shown 
in \cite{frans,vanderMeer2010719} and can be incorporated in our framework without major difficulties.
Besides, inertia effects are also skipped.
In order to trace equilibrium curves we use either a displacement control (for problems without snapbacks)
and the energy-based arc-length control \cite{gutirrez_energy_2004,verhoosel_dissipation-based_2009}.
Interested reader can refer to \cite{Nguyen2013,frans} for the computer implementation aspects of this arc-length
solver. A full Newton-Raphson method was used to solve the discrete equilibrium equations. 
Unless otherwise stated, a geometrically linear formulation is adopted. We use a C++ code \cite{jemjive}
for computations since Matlab is not suitable for this purpose. Whenever possible, validation against 
theoretical solutions are provided.

Four numerical examples are provided including

\begin{itemize}
\item Mixed mode bending test (MMB), 2D simple geometry, implementation verification test;
\item L-shaped specimen, single and multiple delamination, NURBS curved geometry;
\item 3D double cantilever beam, to verify the implementation;
\item Singly curved thick-walled laminate, 3D curved geometry.
\end{itemize}
And in an extra example, we present NURBS parametrization for other commonly used 
composite structures--glare panel with a circular initial delamination, open hole laminate
and doubly curved composite panel.

\subsection{Mixed mode bending test (MMB)}\label{mmb}

\begin{figure}[htbp]
         \centering
         \psfrag{l}[c]{$L$}\psfrag{a0}[c]{$a_0$}\psfrag{h}{$h$}\psfrag{P1}{$P_1,u_1$}\psfrag{P2}{$P_2,u_2$}
         \includegraphics[width=\textwidth]{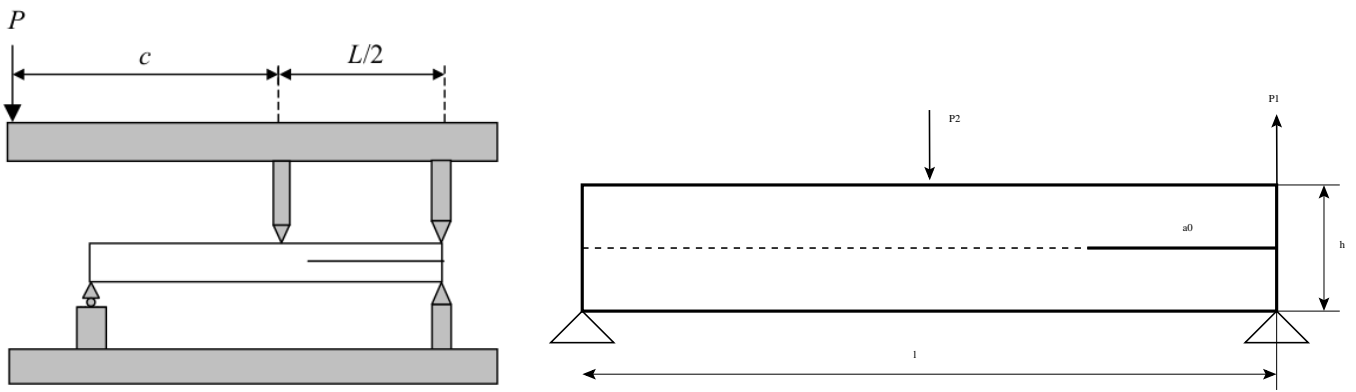}
        \caption{Mixed Mode Bending (MMB): geometry and loading.} \label{fig:mmb}
\end{figure}

Figure \ref{fig:mmb} shows the mixed mode bending test of which the geometry data
are $L=100$ mm, $h=3$ mm; the beam thickness $B$ is equal to 10 mm; the initial crack length is $a_0=20$ mm.  
The plies are modeled with isotropic material to make a fair
comparison with analytical solutions \cite{Mi01071998} which are valid for isotropic materials only. 
The properties for the isotropic material are $E=150$ GPa and $\nu=0.25$. The properties for the cohesive
elements are $G_{Ic}=0.352$ N/mm, $G_{IIc}=1.45$ N/mm and $\tau^0_1=80$ MPa,
$\tau^0_3=60$ MPa. The interface stiffness is $K=10^6$ N/mm$^3$ and $\eta=1.56$.
In order to prevent interpenetration of the two arms, in addition to cohesive elements, 
frictionless contact elements are placed along the initial crack.
The loads applied are $P_1=2Pc/L$ and
$P_2=P(2c+L)/L$, where $L$ is the beam length, $c$ is the lever arm length, and $P$ is the applied
load. From these relationships, it is clear that the applied loads $P_1$ and $P_2$ are proportional
\ie $P_2/P_1=(2c+L)/L$. We choose $c=43.72$ mm so that the mixed-mode ratio $G_I/G_{II}$ is unity.  
The external force vector is therefore $\vm{f}^\text{ext}=\lambda[1, -2.1436]\trans$ (a unit force was assigned to
$P_1$) in which the variable load scale $\lambda$
is solved together with the nodal displacements using the energy based arc-length method
\cite{gutirrez_energy_2004,verhoosel_dissipation-based_2009,Nguyen2013}.

\subsubsection{Geometry and mesh}

For those who are not familiar to B-splines/NURBS, we present how to build the beam geometry using B-splines.
It is obvious that the beam can be exactly represented by a bilinear B-spline surface 
with 4 control points locating at its four corners. The Matlab code for doing this is lines 1--9 in 
Listing \ref{list-mmb}. Next, the B-spline is order elevated to the order that suits the analysis purpose,
see line 10 of the same Listing. The delamination path locates in the midline of the beam \ie $\eta=0.5$
and note that $q=2$, in order to introduce a discontinuity one simply has to insert $0.5$ three ($=q+1$) 
times into knot vector $\mathscr{H}$ (knot vector which is perpendicular to the delamination plane). 
For point load $P_2$ one needs a control point at the location of the
force which corresponding to insert 0.5 three times (equals $p=3$) into knots $\Xi$. Line 13 does exactly that.
In order to differentiate cohesive elements and contact elements (remind that contact elements are put along
the initial crack to prevent interpenetration), a knot $1-a_0/L$ is added to $\Xi$ $p$ times (see line 14).
The final step is to perform a $h$-refinement to refine the mesh and extract element connectivity data for the
interface elements using the code given in Listing \ref{list-1D}.

\begin{snippet1}[caption={Matlab code to build the beam using B-splines}, 
   label={list-mmb},framerule=1pt]
    controlPts          = zeros(4,2,2);
    controlPts(1:2,1,1) = [0;0];       % L=length of beam
    controlPts(1:2,2,1) = [L;0];       % W=height of beam
    controlPts(1:2,1,2) = [0;W];
    controlPts(1:2,2,2) = [L;W];
    controlPts(4,:,:)   = 1;           % weights = 1: B-splines
    uKnot = [0 0 1 1];
    vKnot = [0 0 1 1];
    solid = nrbmak(controlPts,{uKnot vKnot}); % build the initial NURBS object
    solid = nrbdegelev(solid,[2 1]);          % evaluate order to cubic-quadratic surface
    % h-refinement in Y direction to make sure it is C^{-1} along the 
    % delamination path. The multiplicity must be q+1. Also insert knot to have a C^0 
    % at loading point
    solid     = nrbkntins(solid,{[0.5 0.5 0.5] [0.5 0.5 0.5]});
    solid     = nrbkntins(solid,{[1-a0/L 1-a0/L 1-a0/L] []});% insert for intitial crack 
\end{snippet1}

\subsubsection{Analyses with varying basis orders}

\begin{figure}[htbp]
         \centering
  \subfloat[]{\includegraphics[width=0.45\textwidth]{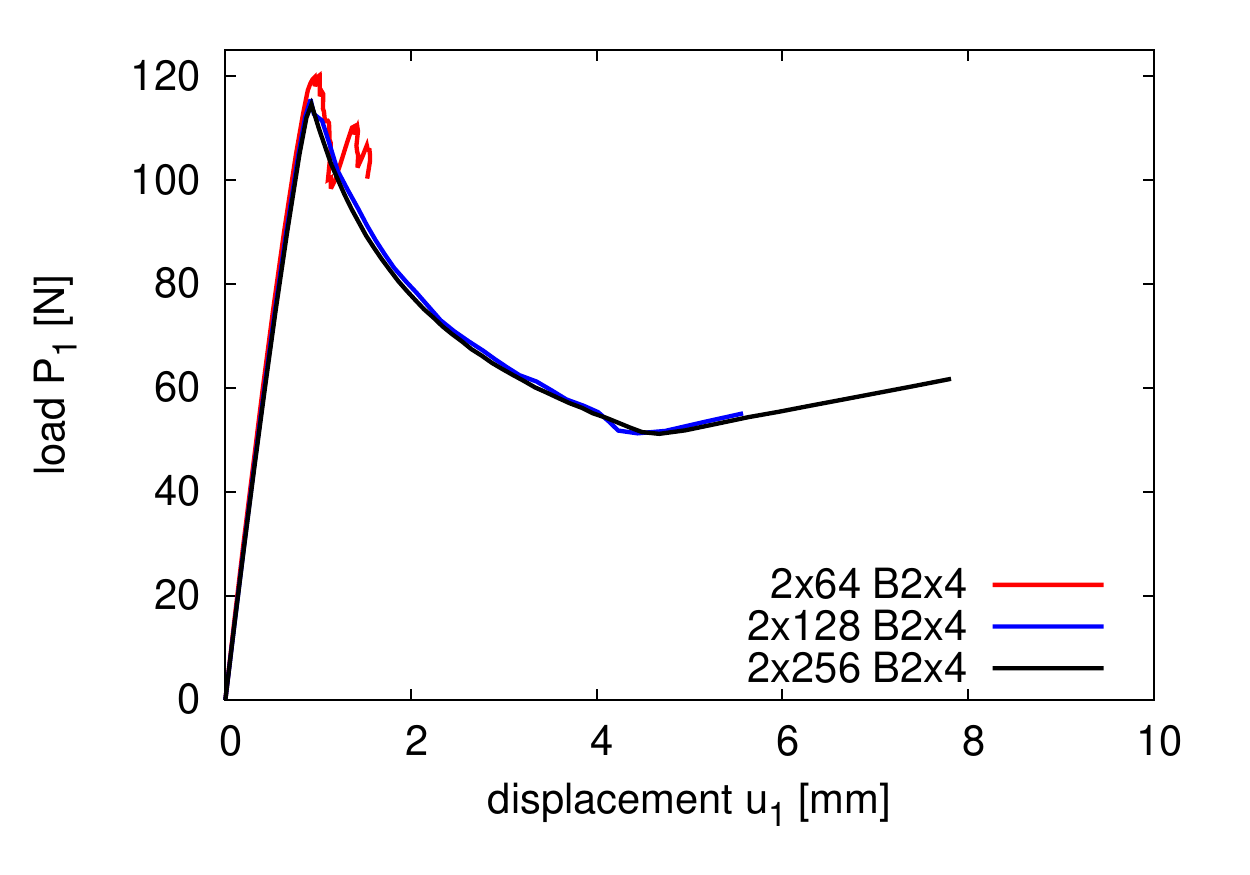}}\;
  \subfloat[]{\includegraphics[width=0.45\textwidth]{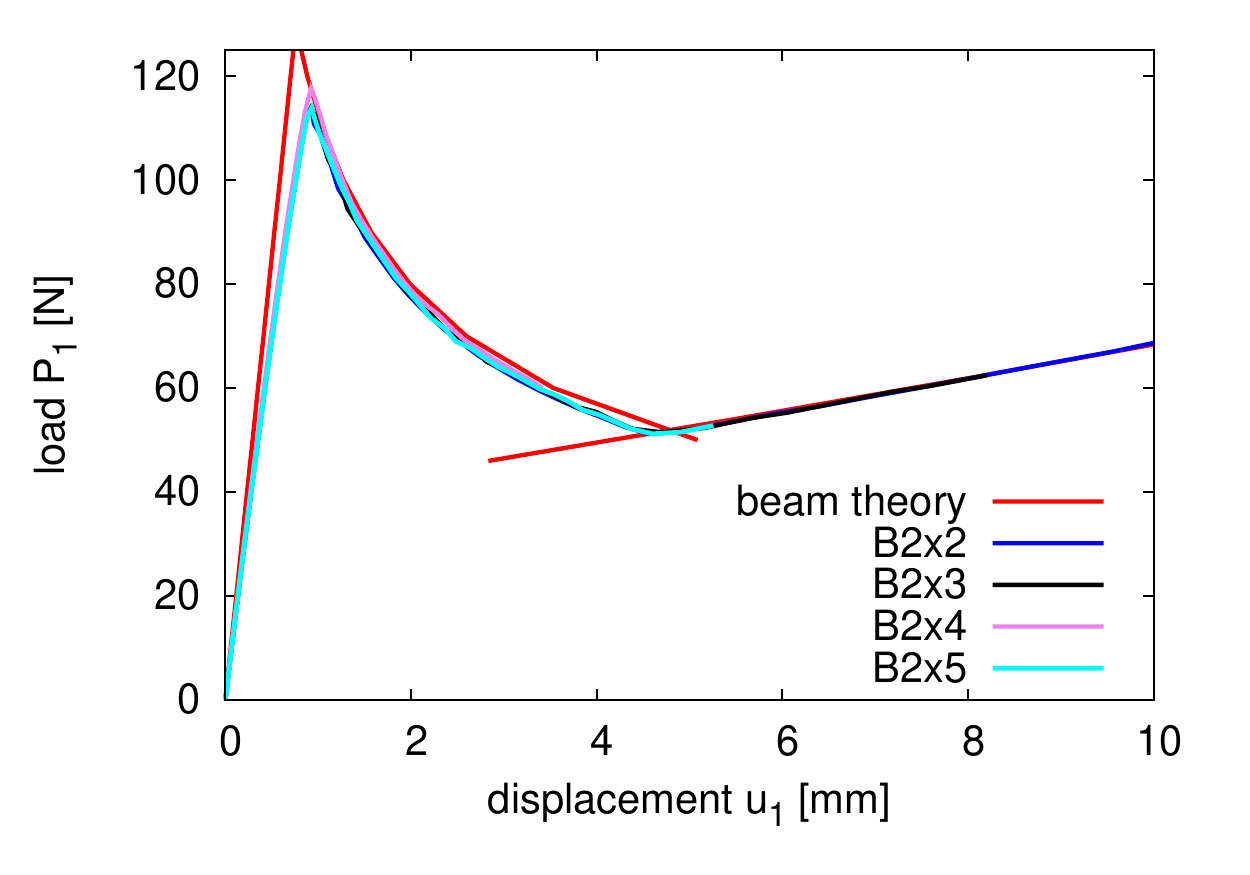}}\;
        \caption{Mixed Mode Bending (MMB): (a) mesh convergence test  and (b)
           varying basis order in the length direction on meshes of $2\times128$ elements.} 
        \label{fig:mmb-lodi}
\end{figure}

We use meshes with two elements along the thickness direction and the basis order along this direction is fixed
to 2 (quadratic basis). The notation $2\times128$ B$2\times3$ indicates a mesh of $2\times128$ elements 
of orders $2\times3$. The order of basis functions along the length direction, $p$, varies from two to five.  
Firstly we perform a mesh convergence test for quartic-quadratic elements and the result is given in 
Fig. (\ref{fig:mmb-lodi}a). 
Mesh $2\times64$ is simply too coarse to accurately capture the cohesive zone and mesh $2\times128$ is
sufficient to get a reasonable result.
Next, the mesh density is fixed at  $2\times128$ and $p$ is varied from 2 to 5,
the result is plotted in Fig. (\ref{fig:mmb-lodi}b). We refer to \cite{Nguyen2013} for a throughout
study on the excellent performance of high order B-splines elements compared to low order Lagrange
finite elements for delamination analyses.

\subsection{L-shaped composite panel with a fillet}

\begin{figure}[htbp]
         \centering
         \psfrag{u}[c]{$u$}
         \psfrag{h}[c]{6.4 mm}
         \psfrag{r}[c]{2.55}
         \psfrag{r1}[c]{2.25}
         \psfrag{x}{$x$}
         \psfrag{y}{$y$}
         \psfrag{layout}{$[0/90/0/90/...]$}
         \includegraphics[width=0.4\textwidth]{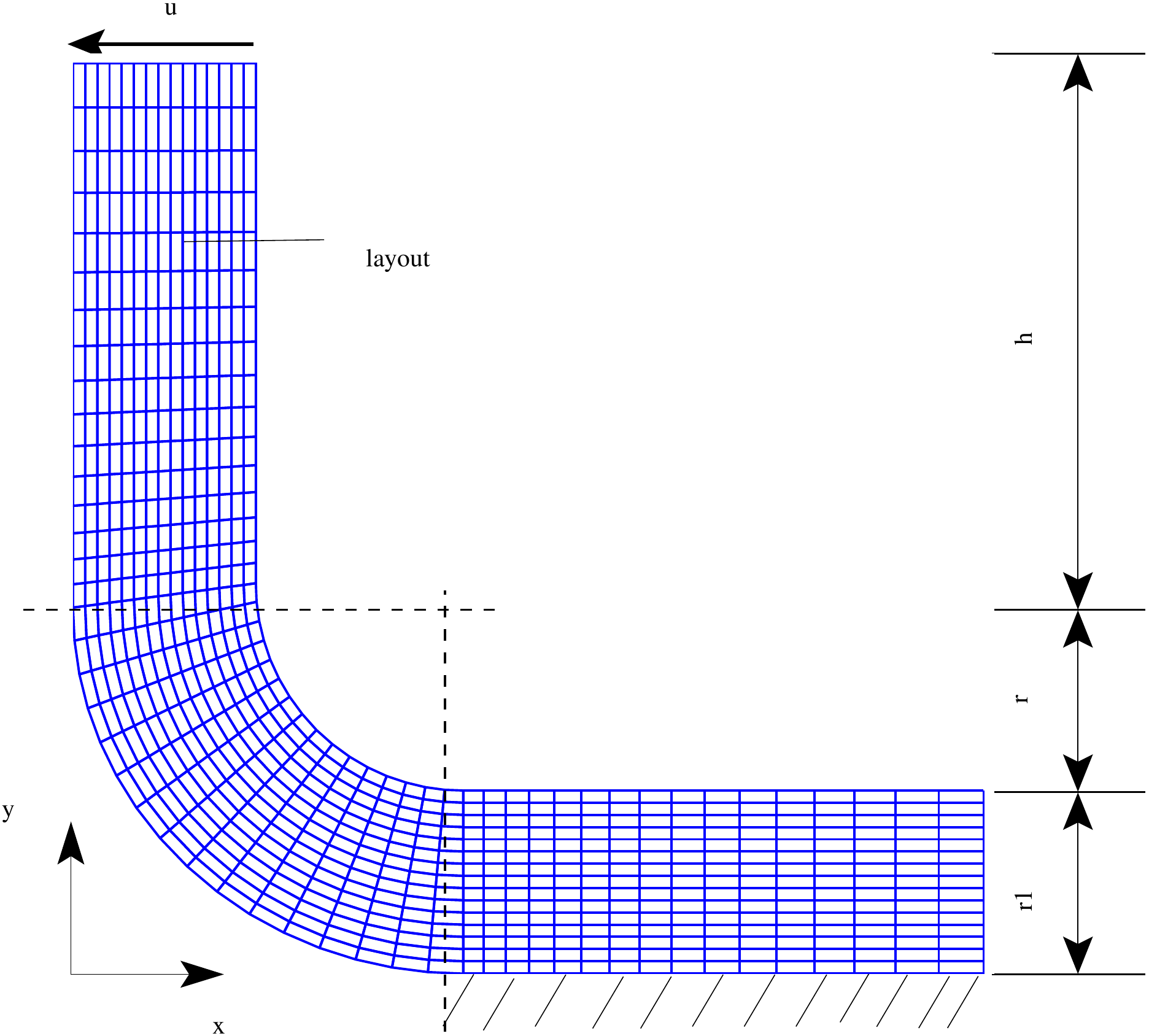}
        \caption{L-shaped specimen: boundary and geometry data. 
          There are 15 plies (0$^\circ$ and 90$^\circ$). 
             The ply orientation is measured with respect to the $x-y$ plane.} 
        \label{fig:lshape-geo}
\end{figure}

For the second example, we analyze the L-shaped composite specimen which was studied in 
\cite{Wimmer20082332,lshape} using Lagrange finite elements. 
The geometry and loading configuration is given in Fig. (\ref{fig:lshape-geo}).
Contrary to the previous example, in this example NURBS surfaces are used to exactly represent the curved
geometry (to be precise circular arcs).
The structure is built up by 15 plies of a unidirectional fiber reinforced carbon/epoxy material. 
The plies are oriented in alternating 0$^\circ$ and 90$^\circ$ orientation, 
where the angle is measured from the $xy$ plane. The inner ply and the outer ply are 
oriented in the 0$^\circ$ direction.
Material constants are given in Table \ref{tab:lshape-material} which are taken from 
\cite{Wimmer20082332,lshape}. A plane strain condition is assumed. For this problem, unless 
otherwise stated, we use bi-quadratic NURBS elements for the continuum and quadratic NURBS
interface elements for the delamination. 

\begin{table}[htbp]
  \centering
  \begin{tabularx}{\textwidth}{XXXXX}
    \toprule
    $E_{11}$  & $E_{22}=E_{33}$  & $G_{12}=G_{13}$  &  $\nu_{12}=\nu_{13}$ & $\nu_{23}$  \\
    \midrule
    139.3 GPa & 9.72 GPa & 5.59 GPa & 0.29  & 0.40   \\ \\
    $G_{Ic}$ & $G_{IIc}$ & $\tau^0_1$ & $\tau^0_3$ & $\mu$ \\
    \midrule
    0.193 N/mm & 0.455 N/mm & 60.0 MPa & 80.0 MPa & 2.0 \\
    \bottomrule
  \end{tabularx}
  \caption{L-shaped specimen: material properties.}
  \label{tab:lshape-material}
\end{table}

\subsubsection{Geometry and mesh}

\begin{figure}[htbp]
         \centering
         \includegraphics[width=0.3\textwidth]{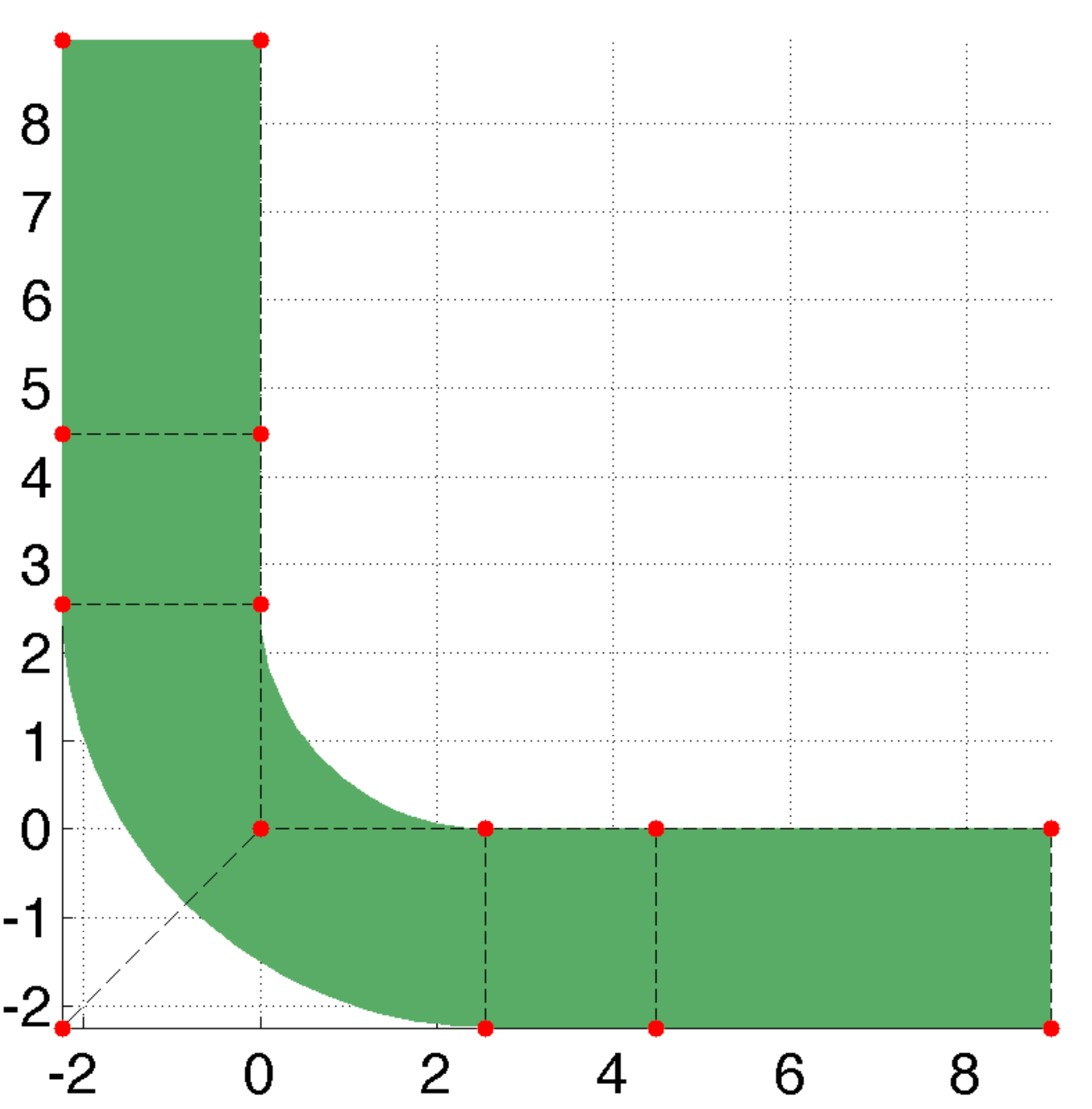}
        \caption{L-shaped specimen: quadratic-linear NURBS geometry with control points 
           (filled circles) and control polygon.} 
        \label{fig:lshape-mesh0}
\end{figure}

The L-shaped geometry can be exactly represented by a quadratic-linear NURBS surface
as shown in Fig. (\ref{fig:lshape-mesh0}) that consists of $7\times2$ control points. 
The Matlab code used to build the NURBS is given in Listing \ref{list:lshape-geo}.
It is easy to vary the number of plies (see line 4 of the same Listing).
Listing \ref{list:lshape-insert} gives code to perform $p$-refinement (to a bi-quadratic NURBS
surface) and knot insertion at ply interfaces (two times) to create $C^{0}$ lines so that
the strain field is discontinuous across the ply interfaces. Next, knot insertion is performed again to generate
discontinuity lines at the desired ply interfaces. Two cases are illustrated in the code--interface
elements locate along the interface between ply 5 and 6 (line 10) and interface elements at every
ply interfaces (line 12-16).

\begin{snippet1}[caption={Matlab code to build NURBS geometry of the L-shaped specimen.}, 
   label={list:lshape-geo},framerule=1pt]
    H  = 6.4;
    R  = 2.55;
    R0 = 2.25;
    no = 15;   % number of plies
    t  = R0/no;% ply thickness
    % initial mesh- quadratic x linear
    uKnot = [0 0 0 1 1 2 2 3 3 3]; uKnot = uKnot/max(uKnot);
    vKnot = [0 0 1 1];
    controlPts      = zeros(4,7,2);
    controlPts(1:2,1,1) = [H+R;0];           controlPts(1:2,1,2) = [H+R;-R0];
    controlPts(1:2,2,1) = [(H+R)/2;0];       controlPts(1:2,2,2) = [(H+R)/2;-R0];
    controlPts(1:2,3,1) = [R;0];             controlPts(1:2,3,2) = [R;-R0];
    controlPts(1:2,4,1) = [0;0];             controlPts(1:2,4,2) = [-R0;-R0];
    controlPts(1:2,5,1) = [0;R];             controlPts(1:2,5,2) = [-R0;R];
    controlPts(1:2,6,1) = [0;(H+R)/2];       controlPts(1:2,6,2) = [-R0;(H+R)/2];
    controlPts(1:2,7,1) = [0;H+R];           controlPts(1:2,7,2) = [-R0;H+R];
    fac = 1/sqrt(2);
    controlPts(4,:) = 1;                       % all weights are units except two points
    controlPts(4,4,1) = controlPts(4,4,2) = fac;
    controlPts(1:2,4,1) = fac*controlPts(1:2,4,1);
    controlPts(1:2,4,2) = fac*controlPts(1:2,4,2);
    solid = nrbmak(controlPts,{uKnot,vKnot});   %% build NURBS object
\end{snippet1}

\begin{snippet1}[caption={Matlab code to generate discontinuity lines.}, 
   label={list:lshape-insert},framerule=1pt]
    solid = nrbdegelev(solid,[0 1]); % evaluate order to bi-quadratic
    % knot insertion to have C^0 lines at ply interfaces
    knots=[];
    for ip=1:no-1
        dd = t*ip/R0;
        knots = [knots dd dd]; % multiplicity = order + 1
    end
    solid     = nrbkntins(solid,{[] knots});   
    % knot insertion to have C^{-1} line at interface between plies 5 and 6
    solid     = nrbkntins(solid,{[] 5*t/R0});  
    % if interface elements are placed at all ply interfaces, then
    knots=[];
    for ip=1:no-1
        dd = t*ip/R0;
        knots = [knots dd]; 
    end
    solid     = nrbkntins(solid,{[] knots});   
    % h-refinement along \xi direction to have a refined model for FEA
\end{snippet1}

\subsubsection{Single delamination with and without initial cracks}

\begin{figure}[htbp]
         \centering
         \includegraphics[width=0.55\textwidth]{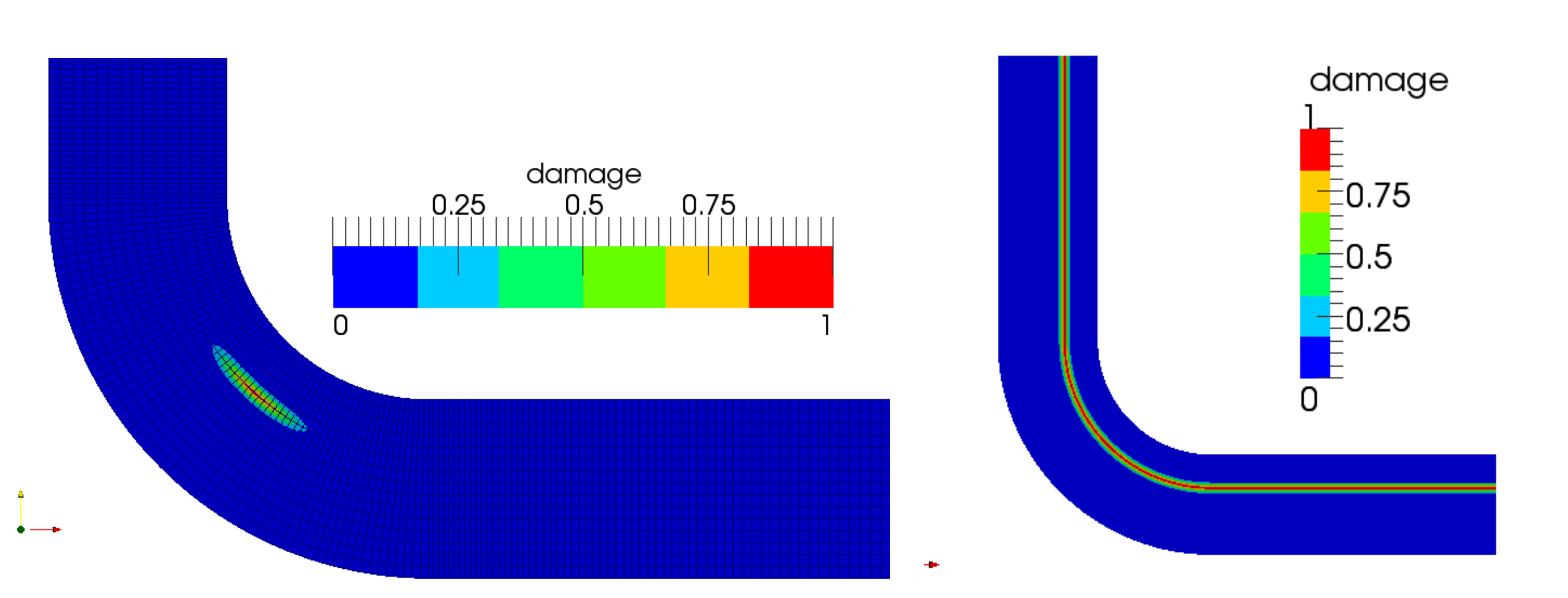}
        \caption{L-shaped specimen: delamination configurations at the peak (left) and when
           the delamination reached the two ends (right).}
        \label{fig:lshape-damage1}
\end{figure}

\begin{figure}[htbp]
         \centering
         \includegraphics[width=0.45\textwidth]{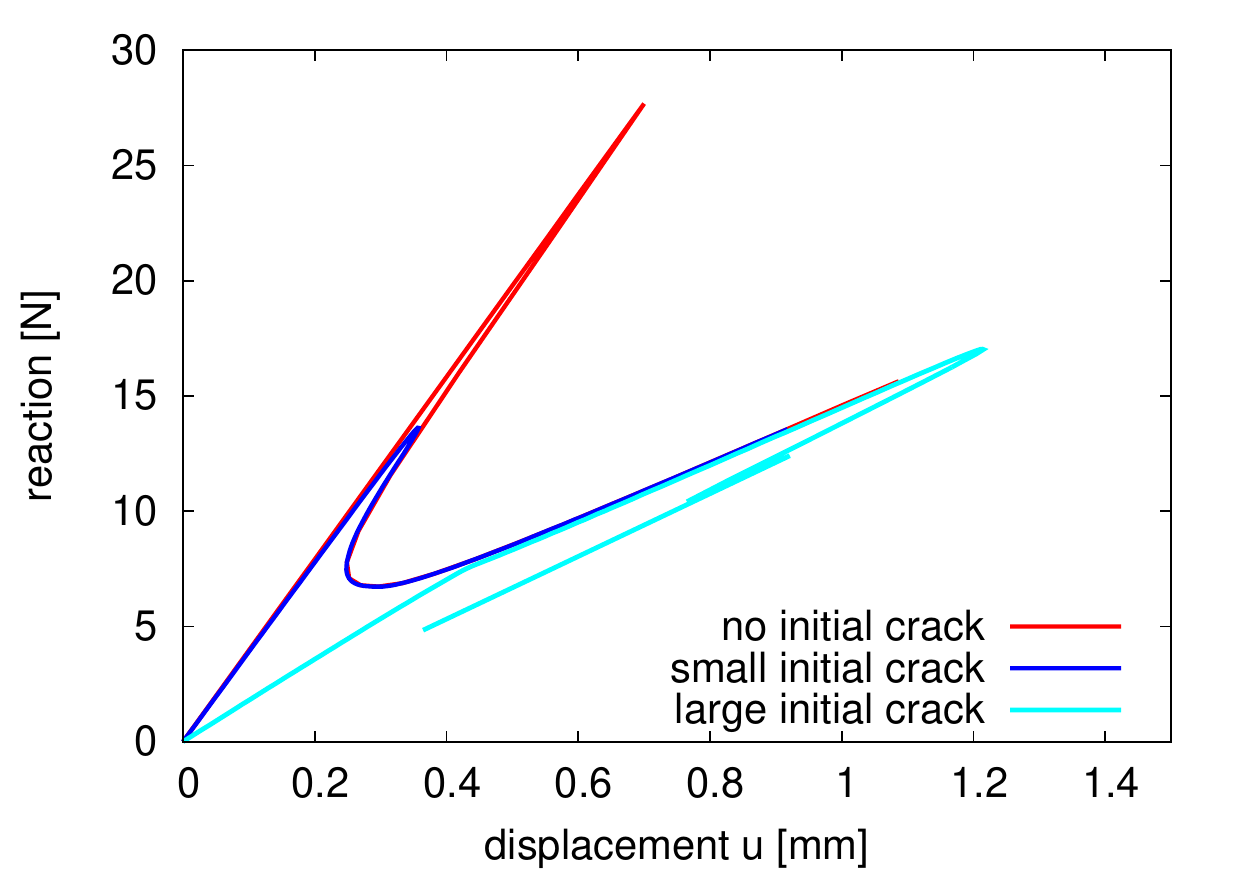}
        \caption{L-shaped specimen with one single delamination between ply 5 and 6:
                without initial cracks, with an small and large initial crack. }
        \label{fig:lshape-lodi}
\end{figure}

Delamination of the interface between ply five and six is first analyzed.
Note that at other ply interfaces, there is no cohesive elements. 
Firstly, the case of no initial cracks is considered. 
One layer of elements is used for each ply. 
The contour plot of damage on the deformed shape is given in Fig. (\ref{fig:lshape-damage1})
and the response in terms of reaction-displacement curve is plotted in Fig. (\ref{fig:lshape-lodi}).
There is a sharp snap-back that corresponds to an unstable delamination growth. After the delamination reaches
a certain size, stable delamination growth is observed as shown by the increasing part of the load-displacement
curve. This is in good agreement with the semi-analytical analysis in \cite{Wimmer20082332}.
The excellent performance of the energy-based arc-length control for responses with sharp snap-backs 
has been demonstrated elsewhere \eg \cite{Nguyen2013,frans}, we therefore do not give an 
discussion on this issue.

\begin{figure}[htbp]
         \centering
         \psfrag{x1}{$\vm{x}_1$}
         \psfrag{x2}{$\vm{x}_2$}
         \psfrag{init}{initial crack}
         \includegraphics[width=0.35\textwidth]{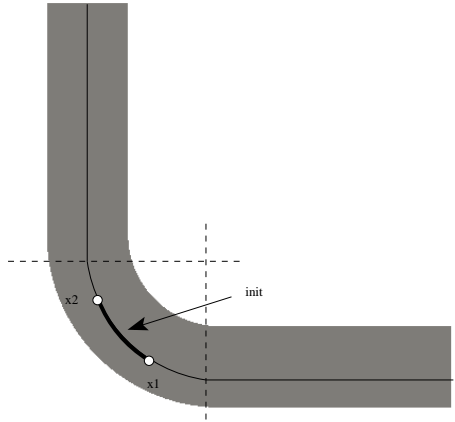}
        \caption{L-shaped specimen with one initial crack.}
        \label{fig:lshape-init-geo}
\end{figure}

Let assume now that there is an initial crack lying on the interface between ply 5 and 6,
see Fig. (\ref{fig:lshape-init-geo}). The initial crack is a part of the NURBS curve that
defines the interface of ply 5 and 6. In this case the geometry modeling procedure is more
involved and follows the steps given in Listing \ref{list:lshape-init-crack}. The extra step
is to perform a point inversion algorithm \cite{piegl_book} to find out the parameters $\xi_1$
and $\xi_2$ that correspond to points $\vm{x}_1$ and $\vm{x}_2$--the tips of the initial crack. 
After that $\xi_1,\xi_2$
are inserted twice (remind that the NURBS basis order along the $\xi$ direction is two).

\begin{figure}[htbp]
         \centering
         \psfrag{a}{(a) small initial crack}
         \psfrag{b}{(b) large initial crack}
         \includegraphics[width=0.5\textwidth]{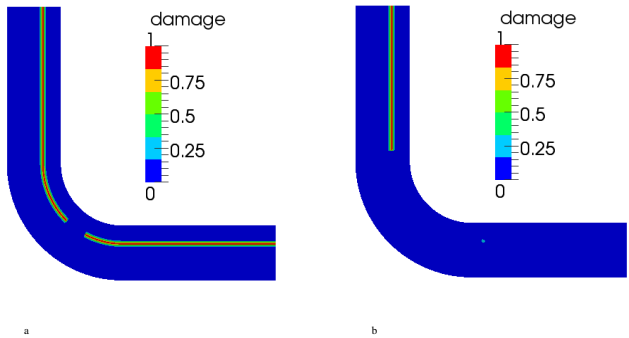}
        \caption{L-shaped specimen with one initial crack: delamination configurations for the case of
        small initial crack (a) and large initial crack (b).}
        \label{fig:lshape-damage2}
\end{figure}

\begin{snippet1}[caption={L-shaped specimen with an initial crack: Matlab code to build the geometry.}, 
   label={list:lshape-init-crack},framerule=1pt]
    % code from Listing 4 to build the NURBS surface
    % code from Listing 5 to create C^0 and C^{-1} lines
    % point inversion to find parametric values xi1 and xi2 that correspond to
    % points x_1 and x_2 defining the initial crack.
    % insert xi1 and xi2 2 times (p=) to have C^0 at x_1 and x_2
    solid     = nrbkntins(solid,{[xi1 xi1 xi2 xi2] []});   
    % h-refinement along \xi direction to have a refined model for FEA
\end{snippet1}

\begin{rmk}
Point inversion for NURBS curves concerns the computation of parameter $\bar{\xi}$ that
corresponds to a point $\bar{\vm{x}}$ such that $N_I(\bar{\xi})\vm{x}_I=\bar{\vm{x}}$ where
$\vm{x}_I$ denote the control points of the curve. Generally, a Newton-Raphson iterative
method is used, we refer to \cite{piegl_book} for details.
\end{rmk}

Two cases, one with a small initial crack and one with a large initial crack are considered.
The delamination of the specimen is given in Fig. (\ref{fig:lshape-damage2}) and the responses
are plotted in Fig. (\ref{fig:lshape-lodi}). For the case of a small initial crack, the response
of the specimen is very similar to the case without any initial cracks, except that the peak load
is smaller.
For the case of a large initial crack, the delamination growth is stable.
This is in good agreement with the work in \cite{Wimmer20082332}.

\subsubsection{Multiple delaminations}

\begin{figure}[htbp]
         \centering
         \includegraphics[width=0.75\textwidth]{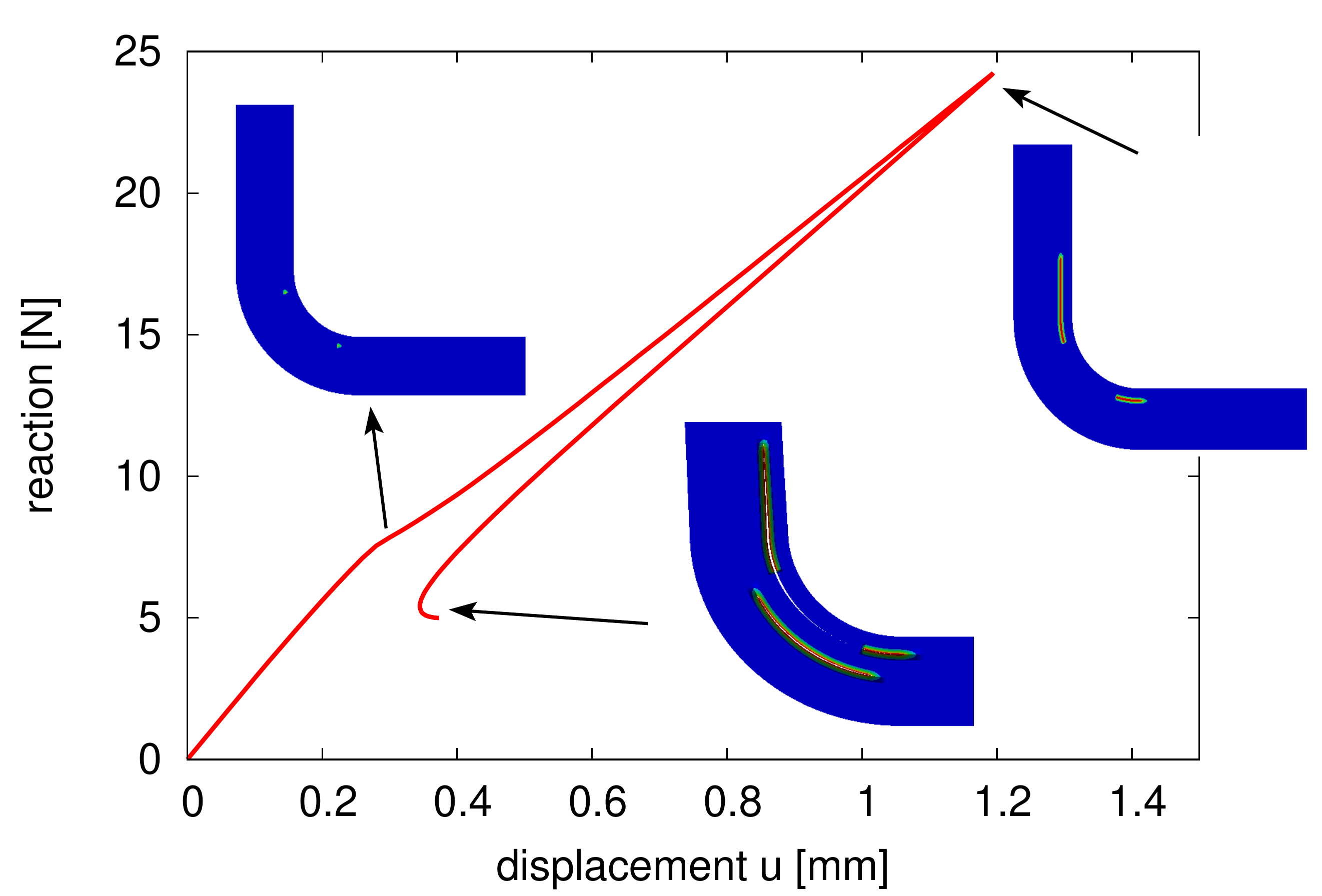}
        \caption{L-shaped specimen with one initial crack and cohesive elements at all ply interfaces: 
           multiple delaminations.}
        \label{fig:lshape-multi-lodi}
\end{figure}

In order to study multiple delaminations, we place cohesive elements along all ply 
interfaces and one large initial crack at the interface
of ply 3 and 4 (we conducted an analysis without any initial crack and found that delamination was
initiated at the interface of ply 3 and 4). The analysis was performed using about 100 load increments
and the computation time was 730s on a Intel Core i7 2.8GHz laptop (29 340 unknowns and 9280 elements).
Figure \ref{fig:lshape-multi-lodi} gives the response of  the specimen. As can be seen, the propagation
of the first delamination (from both tips of the initial crack) is stable and the growth of the second
delamination (between ply 7 and 8) is unstable \footnote{Movies of these analyses can be found at
\url{http://www.frontiersin.org/people/NguyenPhu/94150/video}.}.

\subsection{Three dimensional double cantilever beam}

\begin{figure}[htbp]
         \centering
         \psfrag{p}{$P$}
         \psfrag{le}{length $L$: 100 mm}
         \psfrag{wi}{width $W$: 20mm}
         \psfrag{th}{thickness $b$: 3 mm}
         \psfrag{cr}{initial crack length $a_0$: 30 mm}
         \includegraphics[width=0.45\textwidth]{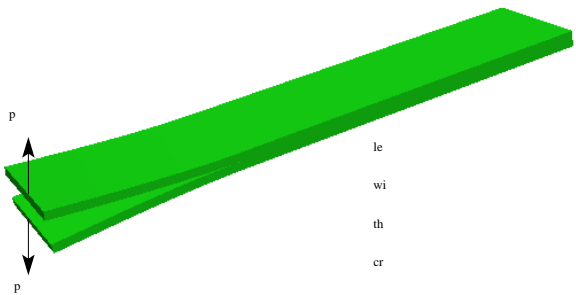}
        \caption{Three dimensional double cantilever beam: geometry and loading data.} 
        \label{fig:dcb-geo}
\end{figure}

As the simplest 3D delamination problem as far as geometry is concerned, we consider the 3D
double cantilever beam (DCB) problem as given in Fig. (\ref{fig:dcb-geo}). This example serves
as a verification test for (1) verifying the implementation of 3D isogeometric interface elements
and (2) validating the automatic generation of 2D isogeometric interface elements.

\subsubsection{Geometry and mesh}

The beam geometry is represented by one single tri-linear NURBS (actually B-splines as the weights
are all units), see Lines 1--5 of Listing \ref{list:dcb-geo}. Order elevation was then performed
to obtain a tri-quadratic solid (line 7) followed by a knot insertion to create the discontinuity 
surface. Finally, $h$-refinement can be applied along any or all directions to have a refined 
model which is analysis suitable. 
Listing \ref{list-2d} is then used to build the element connectivity array for the interface elements.

\begin{snippet1}[caption={Matlab code to build NURBS geometry of the 3D DCB},
   label={list:dcb-geo},framerule=1pt]
    uKnot = [0 0 1 1];
    vKnot = [0 0 1 1];
    wKnot = [0 0 1 1];
    %% build NURBS object
    solid = nrbmak(controlPts,{uKnot vKnot wKnot}); % controlPts = 8 corners of the beam
    % evaluate order 
    solid = nrbdegelev(solid,[1 1 1]); % to tri-quadratic NURBS solid
    % h-refinement in thickness direction to make sure it is C^{-1} along the 
    % delamination path. The multiplicity must be order+1.
    solid     = nrbkntins(solid,{[] [] [0.5 0.5 0.5]});   
    % insert a knot at the tip of the initial crack
    solid     = nrbkntins(solid,{[a0/L] [] []});   
\end{snippet1}

\subsubsection{Analysis results}

\begin{figure}[htbp]
         \centering
         \includegraphics[width=0.6\textwidth]{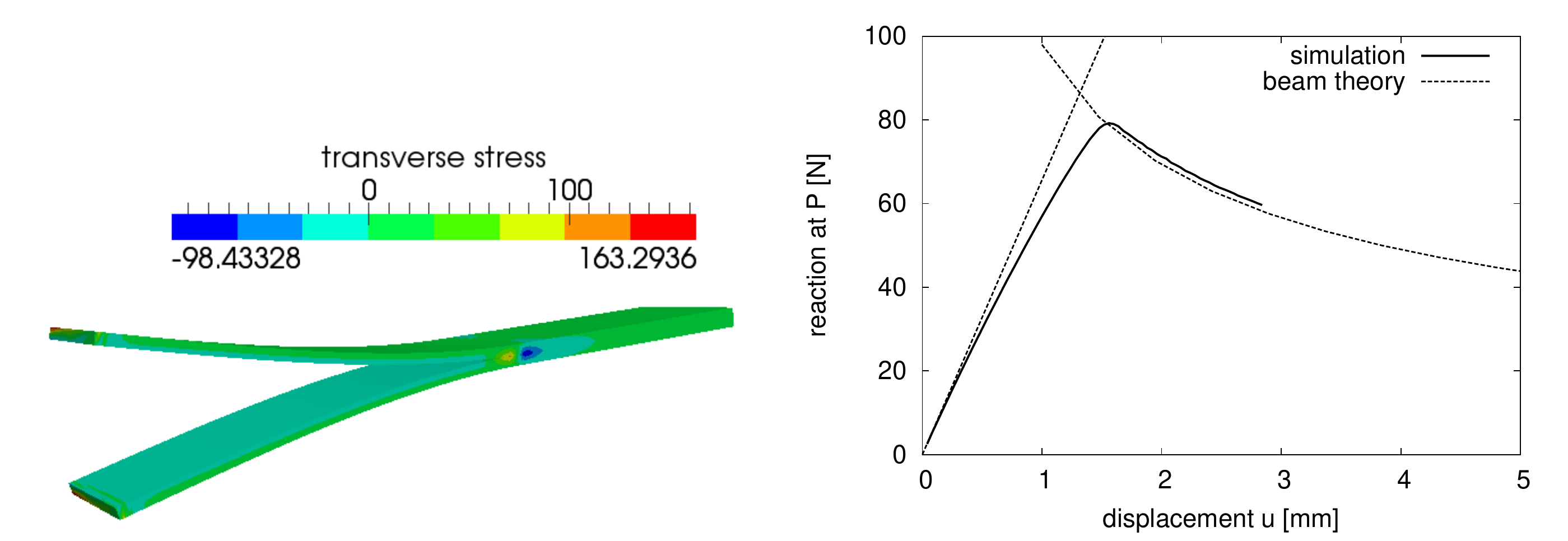}
        \caption{Three dimensional double cantilever beam: contour plot of the transverse stress
           on the deformed shape (magnification factor of 3).} 
        \label{fig:dcb-stress}
\end{figure}

We use an isotropic material with Young modulus $E = 2.1 \times 10^5$ MPa and Poisson ratio $\nu$ = 0.3.
The material constants for the cohesive law are $G_{Ic}$=0.28 N/mm, $\tau^0_1$= 27 MPa, $K=10^7$ N/mm$^3$.
Two layers of elements are placed along the thickness and the width direction.
Figure \ref{fig:dcb-stress} shows the deformed shape and the load-displacement curve including a
comparison with the classical beam theory solution.

\subsubsection{Analysis results with shell elements}

\begin{figure}[htbp]
         \centering
         \psfrag{fix}{fixed displacements}\psfrag{a}{(a)}\psfrag{b}{(b)}\psfrag{c}{(c)}
         \includegraphics[width=0.7\textwidth]{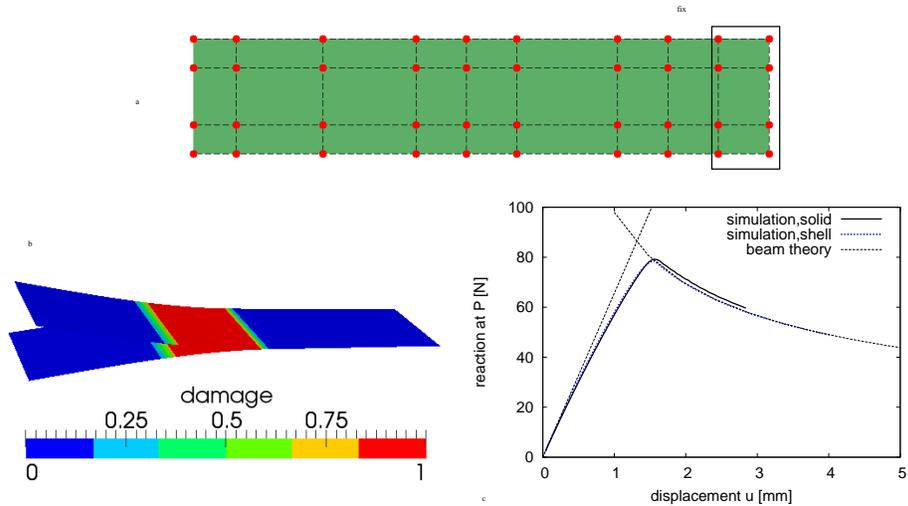}
        \caption{Three dimensional double cantilever beam modeled by shell elements: (a) B-spline mesh
           of one ply with fixed control points and (b) contour plot of damage
           on the deformed shape (magnification factor of 3) and (c) load-displacement curves.} 
        \label{fig:dcb-shell}
\end{figure}

Next, the problem is solved using isogeometric shell elements.
We refer to, for instance, \cite{benson_isogeometric_2010,kiendl_isogeometric_2009,benson_large_2011}
for details on isogeometric shell elements. In this section, we adopt the rotation-free Kirchhoff-Love
thin shell as presented in \cite{kiendl_isogeometric_2009}. Due to its high smoothness, NURBS are suitable
for constructing $C^1$ shell elements without rotation degrees of freedom. In order to fix the rotation
at the right end of the beam, we fix the displacements (all components) of the last two rows of control
points, see Fig. (\ref{fig:dcb-shell}a) and we refer to 
\cite{kiendl_isogeometric_2009,nguyen_iga_review} for details. For each ply is represented by its own
NURBS surface, there is automatically a discontinuity between their interface. Therefore, generation of
interface elements in this context is straightforward.
Each ply is discretized by a mesh of $264\times1$ bi-quadratic elements. The number of nodes/control points
is 1596. Figure \ref{fig:dcb-shell} gives the contour plot of damage and the load-displacement curves.

\subsection{Singly curved thick-walled laminate}

\begin{figure}[htbp]
         \centering
         \psfrag{l}{$l$}
         \psfrag{l}{$l$}
         \psfrag{L}{$L$}
         \psfrag{h}{$h$}
         \psfrag{t}{$t$}
         \psfrag{w}{$w$}
         \includegraphics[width=0.6\textwidth]{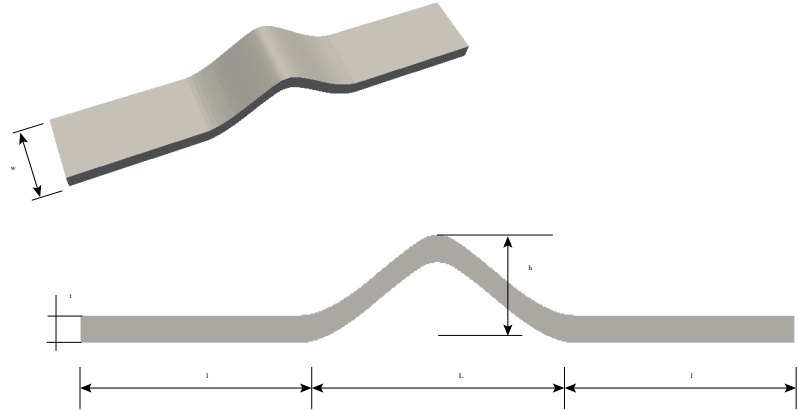}
        \caption{Singly curved thick-walled laminates: geometry configuration. The thickness $t$
        is constant.}
        \label{fig:curved-geo}
\end{figure}

As a 3D example with more complex geometry, we consider a singly curved thick-walled laminate 
which was studied in \cite{Kress2005458,Roos2007327}.
Air-intakes of formula race cars
and strongly curved regions of ship hulls provide examples for such thick-walled curved laminates designs.
The geometry of the sample is given in Fig. (\ref{fig:curved-geo}).
Since the geometry representation of the object of interest is the same in both CAD and FEA
environment, it is very straightforward and fast to get an analysis-suitable model when changes 
are made to the CAD model, for instance changing the thickness $t$. 
This is in sharp contrary to Lagrange finite elements which uses a different geometry representation.
This example also shows how
a trivariate NURBS representation of a curved thick/thin-walled laminate can be built given a
NURBS curve or surface.
For computation, the material constants given in Table \ref{tab:thick-material}
are used of which the material properties of the plies are taken from \cite{Kress2005458}.
The material constants for the interfaces are not provided in \cite{Kress2005458}, the ones used
here are therefore only for computation purposes.

\begin{table}[htbp]
  \centering
  \begin{tabularx}{\textwidth}{XXXXX}
    \toprule
    $E_{11}$  & $E_{22}=E_{33}$  & $G_{12}=G_{13}$  &  $\nu_{12}=\nu_{13}$ & $\nu_{23}$  \\
    \midrule
    110 GPa & 10 GPa & 5.00 GPa & 0.27  & 0.30   \\ \\
    $G_{Ic}$ & $G_{IIc}$ & $\tau^0_1$ & $\tau^0_3$ & $\mu$ \\
    \midrule
    0.075 N/mm & 0.547 N/mm & 80.0 MPa & 90.0 MPa & 1.75 \\
    \bottomrule
  \end{tabularx}
  \caption{Singly curved thick-walled laminate: material properties.}
  \label{tab:thick-material}
\end{table}

\subsubsection{Geometry and mesh}

\begin{figure}[htbp]
         \centering
         \psfrag{base}{base curve}
         \psfrag{o}{offset curve}
         \includegraphics[width=0.6\textwidth]{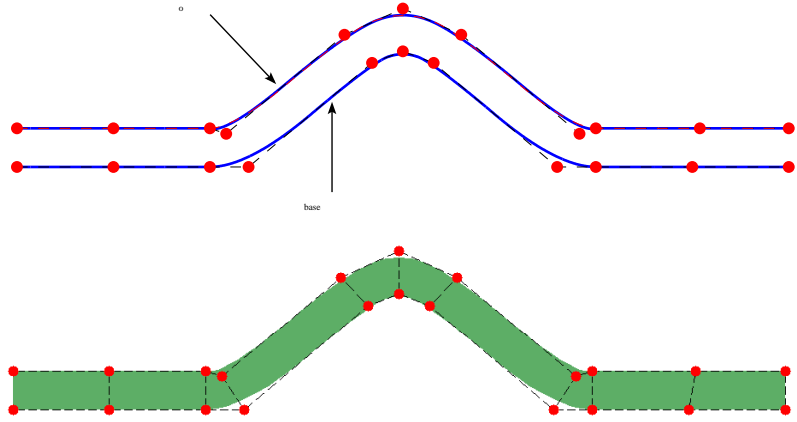}
        \caption{Singly curved thick-walled laminates: building the cross section as a B-spline surface made of
           the base curve and its offset. The red points denote the control points.}
        \label{fig:curved-geom}
\end{figure}

\begin{figure}[htbp]
         \centering
         \psfrag{w}{$w$}
         \includegraphics[width=0.6\textwidth]{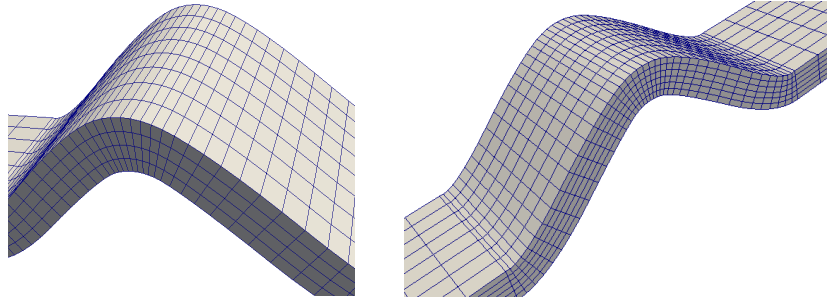}
        \caption{Singly curved thick-walled laminates: NURBS meshes.}
        \label{fig:curved-mesh}
\end{figure}

The geometry of the singly curved thick-walled laminates can be built by first
creating a NURBS curve as shown in Fig. (\ref{fig:curved-geom}). 
Next, an offset of this curve with offset distance $t$ is created using the algorithm
described in \cite{nguyen-offset} which ensures the offset curve has the same parametrization
as its base. This allows a tensor-product surface bounded by these two curves can be constructed.
Having these two curves, a B-spline surface can be constructed.
Knot insertion was then made to build $C^{-1}$ lines at the ply 
interfaces. Finally, the cross section is extruded along the width direction. We refer to
Listing \ref{list:curved-geo} for the Matlab code that produces the geometry. Again, the number of ply
can be easily changed and interface elements can be placed along any ply interface. NURBS meshes
of the sample are given in Fig. (\ref{fig:curved-mesh}).

\begin{snippet1}[caption={Matlab code to build NURBS geometry of the singly curved thick-walled laminates.},
   label={list:curved-geo},framerule=1pt]
    l = 80; L = 100; t = 10; w = 40; h = 30;
    no = 4; % number of layers
    uKnot = [0 0 0 1 1 2 3 4 5  6 6 7 7 7];
    controlPts         = zeros(4,11);
    controlPts(1:2,1)  = [0;0];
    controlPts(1:2,2)  = [0.5*l;0];
    controlPts(1:2,3)  = [l;0];
    controlPts(1:2,4)  = [l+10;0];
    controlPts(1:2,5)  = [l+0.5*L-8;h-3];
    controlPts(1:2,6)  = [l+0.5*L;h];
    controlPts(1:2,7)  = [l+0.5*L+8;h-3];
    controlPts(1:2,8)  = [l+L-10;0];
    controlPts(1:2,9)  = [l+L;0];
    controlPts(1:2,10) = [1.5*l+L;0];
    controlPts(1:2,11) = [2*l+L;0];
    controlPts(4,:,:)  = 1;
    
    curve  = nrbmak(controlPts,uKnot);
    [oCurve,offsetPts]= offsetCurve(curve,t,alpha,beta,eps1,maxIter); % offset curve
    surf   = surfaceFromTwoCurves(curve, oCurve);
    surf   = nrbdegelev(surf,[0 1]); % evaluate order => bi-quadratic
    % h-refinement in Y direction to make sure it is C^{-1} along delamination path
    knots=[];
    for ip=1:no-1
        dd = ip/4;
        knots = [knots dd dd dd];
    end
    surf  = nrbkntins(surf,{[] knots});
    solid = nrbextrude(surf, [0,0,w]);  % make the solid   
\end{snippet1}

\subsubsection{2D analyses}

\begin{figure}[htbp]
         \centering
         \psfrag{w}{$w$}
         \includegraphics[width=0.67\textwidth]{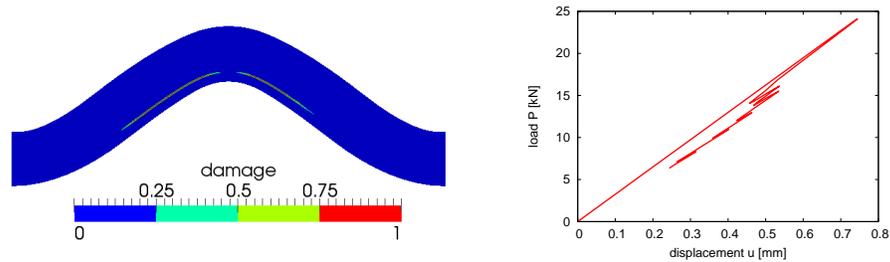}
        \caption{Singly curved thick-walled laminate under compression: delamination pattern (left) and
                load-displacement curve (right).}
        \label{fig:curved-2d-compression}
\end{figure}

Since the straight specimen ends were placed in the clamps of a closed-loop controlled servo-hydraulic 
testing machine \cite{Roos2007327}, 
        in the FE model, the straight ends are not included. The sample is subjected to a compressive
force on the right end and fixed in the left end. The laminate is composed of 45 unidirectional ($0^\circ$) plies of carbon fiber reinforced plastic. The mesh was consisted of $40\times45$ quartic-quadratic NURBS elements
and 1760 quartic interface elements. The number of nodes is 7 020 hence the number of unknowns is 14 040.
The analysis was performed in 121 load increments and the computation time was 1300s 
on a Intel Core i7 2.8GHz laptop. The delamination pattern and the load-displacement is given in 
Fig. (\ref{fig:curved-2d-compression}). Note that no effort was made to compare the obtained result with
the test given in \cite{Roos2007327} since it is beyond the scope of this paper.

\subsubsection{3D analyses}

\subsection{Some other models}

\begin{figure}[htbp]
         \centering
         \psfrag{a}{(a)}
         \psfrag{b}{(b)}
         \psfrag{c}{(c)}
         \psfrag{d}{(d)}
         \includegraphics[width=0.6\textwidth]{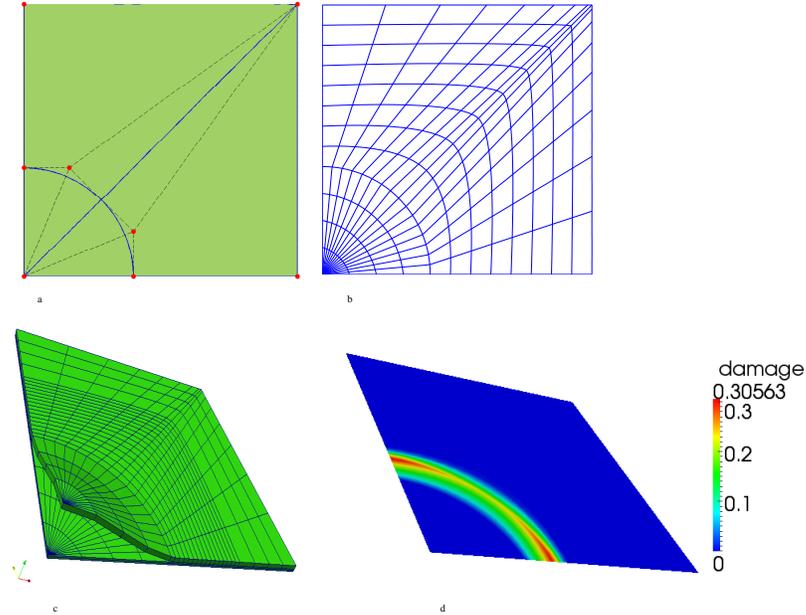}
        \caption{Glare with a circular initial delamination: (a) NURBS surface with control points and
           mesh (4 elements), (b) refined mesh, (c,d) deformed shape and damage plot.}
        \label{fig:circle}
\end{figure}

For completeness, in this section we apply the presented isogeometric framework to other commonly
encountered composite structures. In Fig. (\ref{fig:circle}c), a glare panel with a circular initial
delamination is given (one quarter of the panel is shown due to symmetry). The NURBS representation of the 
panel is given in Fig. (\ref{fig:circle}a) in which the coarsest mesh that consists of $2\times2$ quadratic-linear
NURBS elements can capture the circle geometry and Fig. (\ref{fig:circle}b) shows a refined mesh. Then, interface
elements can be straightforwardly inserted and delamination analyses can be performed  
Fig. (\ref{fig:circle}c,d). It should be emphasized that the chosen NURBS parametrization given in 
Fig. (\ref{fig:circle}a) is not unique and there are singular points at the bottom left and top right corners 
(this, however, does not affect the analysis since no integration points are placed there).

\begin{figure}[htbp]
         \centering
         \psfrag{p1}{patch1}
         \psfrag{p2}{patch2}
         \psfrag{c}{(c)}
         \psfrag{d}{(d)}
         \includegraphics[width=0.6\textwidth]{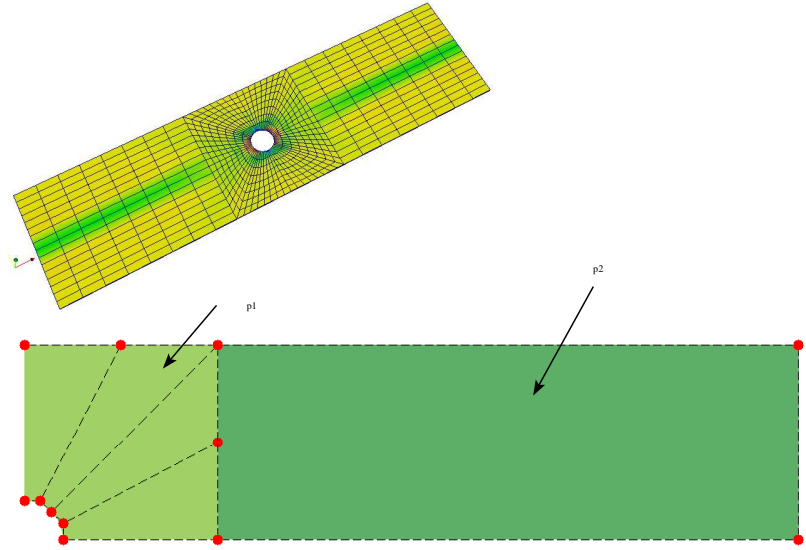}
        \caption{Open hole laminate.}
        \label{fig:hole}
\end{figure}

Next, we present a NURBS mesh for the open hole laminate as shown in Fig. (\ref{fig:hole}).
The whole sample can be represented by six NURBS patches of which four patches are
for the central part. In this figure, we give a parametrization that results in a so-called
compatible multi-patch model. Note that across the patch interface, the basis is only $C^0$.
Interface elements are generated for each patch independently using the presented algorithm.
It should be emphasized that joining NURBS patches of different parametrizations provides 
more flexibility albeit a non-trivial task. T-splines can be used as a remedy, 
see \eg \cite{bazilevs_isogeometric_2010}.

\begin{figure}[htbp]
         \centering
         \psfrag{p1}{patch1}
         \psfrag{p2}{patch2}
         \psfrag{c}{(c)}
         \psfrag{d}{(d)}
         \includegraphics[width=0.75\textwidth]{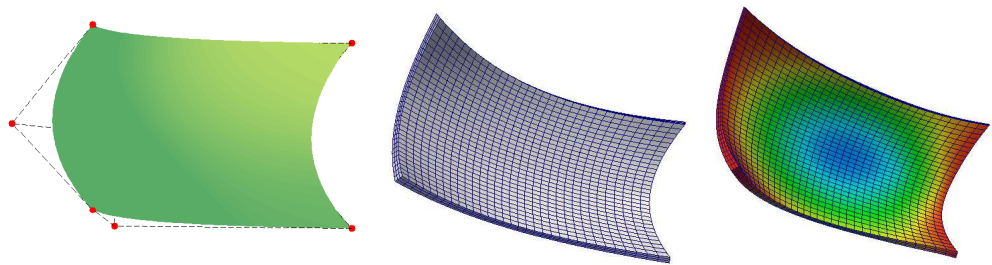}
        \caption{Doubly curved composite laminate: from bi-quadratic NURBS surface to trivariate NURBS solid
           that is suitable for analysis.}
        \label{fig:doubly}
\end{figure}

Finally, treatment of doubly curved composite panels is addresses by an example given in
Fig. (\ref{fig:doubly}). Using a CAD software, the panel is usually a NURBS surface, 
Fig. (\ref{fig:doubly})--left, a trivariate representation is constructed using the ideas
recently reported in \cite{nguyen-offset}, Fig. (\ref{fig:doubly})--middle, 
and FE analyses can be performed, see Fig. (\ref{fig:doubly})--right.

\section{Conclusion}\label{sec:conclusions}

An isogeometric computational framework was presented for modeling delamination of two and three
dimensional composite laminates. By using the isogeometric concept in which the NURBS representation
of the composite laminates is maintained in a finite element environment, it was shown by several examples 
that the time being spent on preparing analysis suitable models for delamination analyses can be dramatically 
reduced. This fact is beneficial to designing composite laminates in which various geometry parameters need to be
varied. The pre-processing algorithms were explained in details and implemented in MIGFEM--an open source code which is available at \url{https://sourceforge.net/projects/cmcodes/}. From the analysis perspective, 
the ultra smooth high order NURBS basis functions are able to produce highly accurate stress fields which is
very important in fracture modeling. The consequence is that relatively coarse meshes (compared to meshes of
lower order elements) can be adopted and thus the computational expense is reduced.

Although an elaborated isogeometric computational framework was presented for modeling delaminated 
composites and several geometry models were addressed, there are certainly a certain number of geometries
that has not been treated. One example is three dimensional curved composite panels with cutouts. Possibilities
for these problems are trimmed NURBS or T-splines for conforming mesh methods and immersed
boundary methods for non-conforming techniques. 

\section*{Acknowledgements}

The authors would like to acknowledge the partial financial support of the
Framework Programme 7 Initial Training Network Funding under grant number
289361 ``Integrating Numerical Simulation and Geometric Design Technology".
St\'{e}phane Bordas also thanks partial funding for his time provided by
1) the EPSRC under grant EP/G042705/1 Increased Reliability for Industrially
Relevant Automatic Crack Growth Simulation with the eXtended Finite Element
Method and 2) the European Research Council Starting Independent Research
Grant (ERC Stg grant agreement No. 279578) entitled ``Towards real time multiscale simulation of cutting in
non-linear materials with applications to surgical simulation and computer
guided surgery''.
The authors would like to express the gratitude towards Drs. Erik Jan Lingen
and Martijn Stroeven at the Dynaflow Research Group, Houtsingel 95, 2719 EB Zoetermeer, The Netherlands 
for providing us the numerical toolkit jem/jive.

\bibliography{isogeometric}
\bibliographystyle{unsrt}


\end{document}

%% file: pre_commands.tex
\makeatletter
\renewcommand\fs@ruled{\def\@fs@cfont{\bfseries}\let\@fs@capt\floatc@ruled
\def\@fs@pre{\hrule height 1.2pt depth0pt \kern2pt}%
\def\@fs@post{\kern2pt\hrule height 1.2pt depth0pt \kern2pt \relax}%
\def\@fs@mid{\kern2pt\hrule\kern2pt}%
\let\@fs@iftopcapt\iftrue}
\makeatother

\newcommand{\cref}[1]{Chapter~\ref{chapt:#1}}

\newcommand{\D}{\displaystyle}

\usepackage{xspace}
\newcommand{\eg}{e.g.,\xspace}

\newcommand{\ie}{i.e.,\xspace}

\newcommand{\etc}{etc.\@\xspace}

\newcommand{\cmatrixb}{\left\{ \begin{matrix}}
\newcommand{\cmatrixe}{\end{matrix} \right\}}

\newcommand{\vm}[1]{\mathbf{#1}}

\newcommand{\bsym}[1]{\boldsymbol{#1}}

\newcommand{\trans}{^\mathrm{T}}

\newcommand{\bc}{\begin{center}}
\newcommand{\ec}{\end{center}}
\newcommand{\bitem}{\begin{itemize}}
\newcommand{\eitem}{\end{itemize}}

\newcommand{\ljump}{\lbrack \! \lbrack } 
\newcommand{\rjump}{\rbrack \! \rbrack } 
\newcommand{\jump}[1]{\ljump {#1} \rjump} 

\newcommand{\beq}{\begin{equation}}
\newcommand{\eeq}{\end{equation}}
\newcommand{\beqa}{\begin{eqnarray}}
\newcommand{\eeqa}{\end{eqnarray}}

\newcommand{\invisible}[1]{}

\newcommand{\bv}{\begin{verbatim}}
\newcommand{\V}{\verb}                  % EX: \V=-d{#@~}= Expr must

\newcommand{\testpix}[1]{\fbox{\begin{minipage}[c]{\textwidth}
                      #1 \end{minipage} }}

\newcommand{\putpstex}[1]{\includegraphics{#1.pstex_t}}

%% file: delamination2.bbl
\begin{thebibliography}{10}

\bibitem{hughes_isogeometric_2005}
{T.J.R.} Hughes, {J.A.} Cottrell, and Y.~Bazilevs.
\newblock Isogeometric analysis: {CAD}, finite elements, {NURBS}, exact
  geometry and mesh refinement.
\newblock {\em Computer Methods in Applied Mechanics and Engineering},
  194(39-41):4135--4195, 2005.

\bibitem{cottrel_book_2009}
J.~A. Cottrell, T.~J.R. Hughes, and Y.~Bazilevs.
\newblock {\em Isogeometric Analysis: Toward Integration of {CAD and FEA}}.
\newblock Wiley, 2009.

\bibitem{NME:NME292}
P.~Kagan, A.~Fischer, and P.~Z. Bar-Yoseph.
\newblock New {B-Spline Finite Element} approach for geometrical design and
  mechanical analysis.
\newblock {\em International Journal for Numerical Methods in Engineering},
  41(3):435--458, 1998.

\bibitem{Cirak_2000}
F.~Cirak, M.~Ortiz, and P.~Schr{\"o}der.
\newblock Subdivision surfaces: a new paradigm for thin-shell finite-element
  analysis.
\newblock {\em International Journal for Numerical Methods in Engineering},
  47(12):2039--2072, 2000.

\bibitem{temizer_contact_2011}
\.{I}. Temizer, P.~Wriggers, and {T.J.R.} Hughes.
\newblock Contact treatment in isogeometric analysis with {NURBS}.
\newblock {\em Computer Methods in Applied Mechanics and Engineering},
  200(9-12):1100--1112, 2011.

\bibitem{jia_isogeometric_2011}
L.~Jia.
\newblock Isogeometric contact analysis: Geometric basis and formulation for
  frictionless contact.
\newblock {\em Computer Methods in Applied Mechanics and Engineering},
  200(5-8):726--741, 2011.

\bibitem{temizer_three-dimensional_2012}
\.{I}. Temizer, P.~Wriggers, and {T.J.R.} Hughes.
\newblock {Three-Dimensional} {Mortar-Based} frictional contact treatment in
  isogeometric analysis with {NURBS}.
\newblock {\em Computer Methods in Applied Mechanics and Engineering},
  209--212:115--128, 2012.

\bibitem{de_lorenzis_large_2011}
L.~De~Lorenzis, \.{I}. Temizer, P.~Wriggers, and G.~Zavarise.
\newblock A large deformation frictional contact formulation using
  {NURBS-bases} isogeometric analysis.
\newblock {\em International Journal for Numerical Methods in Engineering},
  87(13):1278--1300, 2011.

\bibitem{Matzen201327}
M.E. Matzen, T.~Cichosz, and M.~Bischoff.
\newblock A point to segment contact formulation for isogeometric, {NURBS}
  based finite elements.
\newblock {\em Computer Methods in Applied Mechanics and Engineering}, 255:27
  -- 39, 2013.

\bibitem{wall_isogeometric_2008}
W.~A. Wall, M.~A. Frenzel, and C.~Cyron.
\newblock Isogeometric structural shape optimization.
\newblock {\em Computer Methods in Applied Mechanics and Engineering},
  197(33-40):2976--2988, 2008.

\bibitem{manh_isogeometric_2011}
N.~D. Manh, A.~Evgrafov, A.~R. Gersborg, and J.~Gravesen.
\newblock Isogeometric shape optimization of vibrating membranes.
\newblock {\em Computer Methods in Applied Mechanics and Engineering},
  200(13-16):1343--1353, 2011.

\bibitem{qian_isogeometric_2011}
X.~Qian and O.~Sigmund.
\newblock Isogeometric shape optimization of photonic crystals via {Coons}
  patches.
\newblock {\em Computer Methods in Applied Mechanics and Engineering},
  200(25-28):2237--2255, 2011.

\bibitem{xiaoping_full_2010}
X.~Qian.
\newblock Full analytical sensitivities in {NURBS} based isogeometric shape
  optimization.
\newblock {\em Computer Methods in Applied Mechanics and Engineering},
  199(29-32):2059--2071, 2010.

\bibitem{benson_isogeometric_2010}
{D.J.} Benson, Y.~Bazilevs, {M.C.} Hsu, and {T.J.R.} Hughes.
\newblock Isogeometric shell analysis: The {Reissner--Mindlin} shell.
\newblock {\em Computer Methods in Applied Mechanics and Engineering},
  199(5-8):276--289, 2010.

\bibitem{kiendl_isogeometric_2009}
J.~Kiendl, {K.-U.} Bletzinger, J.~Linhard, and R.~W{\"u}chner.
\newblock Isogeometric shell analysis with {Kirchhoff-Love} elements.
\newblock {\em Computer Methods in Applied Mechanics and Engineering},
  198(49-52):3902--3914, 2009.

\bibitem{benson_large_2011}
{D.J.} Benson, Y.~Bazilevs, {M.-C.} Hsu, and {T.J.R.} Hughes.
\newblock A large deformation, rotation-free, isogeometric shell.
\newblock {\em Computer Methods in Applied Mechanics and Engineering},
  200(13-16):1367--1378, 2011.

\bibitem{beirao_da_veiga_isogeometric_2012}
L.~{Beir\~{a}o da Veiga}, A.~Buffa, C.~Lovadina, M.~Martinelli, and
  G.~Sangalli.
\newblock An isogeometric method for the {Reissner-Mindlin} plate bending
  problem.
\newblock {\em Computer Methods in Applied Mechanics and Engineering},
  209--212:45--53, 2012.

\bibitem{uhm_tspline_2009}
{T. K.} Uhm and {S. K.} Youn.
\newblock T-spline finite element method for the analysis of shell structures.
\newblock {\em International Journal for Numerical Methods in Engineering},
  80(4):507--536, 2009.

\bibitem{Echter2013170}
R.~Echter, B.~Oesterle, and M.~Bischoff.
\newblock A hierarchic family of isogeometric shell finite elements.
\newblock {\em Computer Methods in Applied Mechanics and Engineering}, 254:170
  -- 180, 2013.

\bibitem{Benson2013133}
D.J. Benson, S.~Hartmann, Y.~Bazilevs, M.-C. Hsu, and T.J.R. Hughes.
\newblock Blended isogeometric shells.
\newblock {\em Computer Methods in Applied Mechanics and Engineering}, 255:133
  -- 146, 2013.

\bibitem{cottrell_isogeometric_2006}
{J.A.} Cottrell, A.~Reali, Y.~Bazilevs, and {T.J.R.} Hughes.
\newblock Isogeometric analysis of structural vibrations.
\newblock {\em Computer Methods in Applied Mechanics and Engineering},
  195(41-43):5257--5296, 2006.

\bibitem{Hughes20084104}
T.J.R. Hughes, A.~Reali, and G.~Sangalli.
\newblock Duality and unified analysis of discrete approximations in structural
  dynamics and wave propagation: Comparison of p-method finite elements with
  k-method {NURBS}.
\newblock {\em Computer Methods in Applied Mechanics and Engineering},
  197(49--50):4104 -- 4124, 2008.

\bibitem{NME:NME4282}
C.~H. Thai, H.~Nguyen-Xuan, N.~Nguyen-Thanh, T-H. Le, T.~Nguyen-Thoi, and
  T.~Rabczuk.
\newblock Static, free vibration, and buckling analysis of laminated composite
  {Reissner-Mindlin} plates using {NURBS-based} isogeometric approach.
\newblock {\em International Journal for Numerical Methods in Engineering},
  91(6), 2012.

\bibitem{Wang2013}
D.~Wang, W.~Liu, and H.~Zhang.
\newblock Novel higher order mass matrices for isogeometric structural
  vibration analysis.
\newblock {\em Computer Methods in Applied Mechanics and Engineering},
  pages~--, 2013.

\bibitem{gomez_isogeometric_2010}
H.~Gomez, {T.J.R.} Hughes, X.~Nogueira, and V.~M. Calo.
\newblock Isogeometric analysis of the isothermal {Navier-Stokes-Korteweg}
  equations.
\newblock {\em Computer Methods in Applied Mechanics and Engineering},
  199(25-28):1828--1840, 2010.

\bibitem{nielsen_discretizations_2011}
P.~N. Nielsen, A.~R. Gersborg, J.~Gravesen, and N.~L. Pedersen.
\newblock Discretizations in isogeometric analysis of {Navier-Stokes} flow.
\newblock {\em Computer Methods in Applied Mechanics and Engineering},
  200(45-46):3242--3253, 2011.

\bibitem{Bazilevs:2010:LES:1749635.1750210}
Y.~Bazilevs and I.~Akkerman.
\newblock Large eddy simulation of turbulent {Taylor-Couette} flow using
  isogeometric analysis and the residual-based variational multiscale method.
\newblock {\em Journal of Computational Physics}, 229(9):3402--3414, 2010.

\bibitem{bazilevs_isogeometric_2008}
Y.~Bazilevs, V.~M. Calo, T.~J.~R. Hughes, and Y.~Zhang.
\newblock Isogeometric fluid-structure interaction: theory, algorithms, and
  computations.
\newblock {\em Computational Mechanics}, 43:3--37, 2008.

\bibitem{bazilevs_patient-specific_2009}
Y.~Bazilevs, {J.R.} Gohean, {T.J.R.} Hughes, {R.D.} Moser, and Y.~Zhang.
\newblock Patient-specific isogeometric fluid-structure interaction analysis of
  thoracic aortic blood flow due to implantation of the {Jarvik} 2000 left
  ventricular assist device.
\newblock {\em Computer Methods in Applied Mechanics and Engineering},
  198(45-46):3534--3550, 2009.

\bibitem{gomez_isogeometric_2008}
H.~G{\'o}mez, V.~M. Calo, Y.~Bazilevs, and {T.J.R.} Hughes.
\newblock Isogeometric analysis of the {Cahn-Hilliard} phase-field model.
\newblock {\em Computer Methods in Applied Mechanics and Engineering},
  197(49-50):4333--4352, 2008.

\bibitem{verhoosel_isogeometric_2011-1}
C.~V. Verhoosel, M.~A. Scott, T.~J.~R. Hughes, and R.~{de Borst}.
\newblock An isogeometric analysis approach to gradient damage models.
\newblock {\em International Journal for Numerical Methods in Engineering},
  86(1):115--134, 2011.

\bibitem{fischer_isogeometric_2010}
P.~Fischer, M.~Klassen, J.~Mergheim, P.~Steinmann, and R.~M{\"u}ller.
\newblock Isogeometric analysis of {2D} gradient elasticity.
\newblock {\em Computational Mechanics}, 47:325--334, 2010.

\bibitem{Masud2012112}
A.~Masud and R.~Kannan.
\newblock B-splines and {NURBS} based finite element methods for {Kohn-Sham}
  equations.
\newblock {\em Computer Methods in Applied Mechanics and Engineering},
  241-244:112 -- 127, 2012.

\bibitem{nguyen_iga_review}
V.~P. Nguyen, R.~Simpson, S.P.A. Bordas, and T.~Rabczuk.
\newblock Isogeometric analysis: {An} overview and computer implementation
  aspects.
\newblock {\em Advances in Engineering Softwares}, pages~--, 2013.
\newblock submitted.

\bibitem{de_luycker_xfem_2011}
E.~De~Luycker, D.~J. Benson, T.~Belytschko, Y.~Bazilevs, and M.~C. Hsu.
\newblock {X-FEM} in isogeometric analysis for linear fracture mechanics.
\newblock {\em International Journal for Numerical Methods in Engineering},
  87(6):541--565, 2011.

\bibitem{ghorashi_extended_2012}
S.~S. Ghorashi, N.~Valizadeh, and S.~Mohammadi.
\newblock Extended isogeometric analysis for simulation of stationary and
  propagating cracks.
\newblock {\em International Journal for Numerical Methods in Engineering},
  89(9):1069--1101, 2012.

\bibitem{Tambat20121}
A.~Tambat and G.~Subbarayan.
\newblock Isogeometric enriched field approximations.
\newblock {\em Computer Methods in Applied Mechanics and Engineering},
  245--246:1 -- 21, 2012.

\bibitem{Borden201277}
M.~J. Borden, C.~V. Verhoosel, M.~A. Scott, T.~J.R. Hughes, and C.~M. Landis.
\newblock A phase-field description of dynamic brittle fracture.
\newblock {\em Computer Methods in Applied Mechanics and Engineering},
  217--220:77 -- 95, 2012.

\bibitem{verhoosel_isogeometric_2011}
C.~V. Verhoosel, M.~A. Scott, R.~{de Borst}, and T.~J.~R. Hughes.
\newblock An isogeometric approach to cohesive zone modeling.
\newblock {\em International Journal for Numerical Methods in Engineering},
  87(1‐5):336--360, 2011.

\bibitem{Allix199561}
O.~Allix, P.~Ladev\`{e}ze, and A.~Corigliano.
\newblock Damage analysis of interlaminar fracture specimens.
\newblock {\em Composite Structures}, 31(1):61 -- 74, 1995.

\bibitem{Schellekens19931239}
J.C.J. Schellekens and R.~De Borst.
\newblock A non-linear finite element approach for the analysis of mode-{I}
  free edge delamination in composites.
\newblock {\em International Journal of Solids and Structures}, 30(9):1239 --
  1253, 1993.

\bibitem{Crisfield}
G.~Alfano and M.~A. Crisfield.
\newblock Finite element interface models for the delamination analysis of
  laminated composites: mechanical and computational issues.
\newblock {\em International Journal for Numerical Methods in Engineering},
  50(7):1701--1736, 2001.

\bibitem{Krueger200125}
R~Krueger and T.K O'Brien.
\newblock A {shell/3D} modeling technique for the analysis of delaminated
  composite laminates.
\newblock {\em Composites Part A: Applied Science and Manufacturing}, 32(1):25
  -- 44, 2001.

\bibitem{rybicki1977}
E.~F. Rybicki and M.~F. Kanninen.
\newblock A finite element calculation of stress intensity factors by a
  modified crack closure integral.
\newblock {\em Engineering Fracture Mechanics}, 9:931--938, 1977.

\bibitem{kruger2002}
R.~Krueger.
\newblock The virtual crack closure technique: {History}, approach and
  applications.
\newblock Technical report, NASA. NASA/CR-2002-211628, ICASE Report No.
  2002-10, 2002.

\bibitem{Guiamatsia20092640}
I.~Guiamatsia, B.G. Falzon, G.A.O. Davies, and L.~Iannucci.
\newblock Element-free galerkin modelling of composite damage.
\newblock {\em Composites Science and Technology}, 69(15--16):2640 -- 2648,
  2009.

\bibitem{Samimi20112202}
M.~Samimi, J.A.W. van Dommelen, and M.G.D. Geers.
\newblock A self-adaptive finite element approach for simulation of mixed-mode
  delamination using cohesive zone models.
\newblock {\em Engineering Fracture Mechanics}, 78(10):2202 -- 2219, 2011.

\bibitem{Guiamatsia20092616}
I.~Guiamatsia, J.K. Ankersen, G.A.O. Davies, and L.~Iannucci.
\newblock Decohesion finite element with enriched basis functions for
  delamination.
\newblock {\em Composites Science and Technology}, 69(15â€“16):2616 --
  2624, 2009.

\bibitem{mos_finite_1999}
N.~Mo\"{e}s, J.~Dolbow, and T.~Belytschko.
\newblock A finite element method for crack growth without remeshing.
\newblock {\em International Journal for Numerical Methods in Engineering},
  46(1):131--150, 1999.

\bibitem{NME:NME907}
J.~J.~C. Remmers, G.~N. Wells, and R.~Borst.
\newblock A solid-like shell element allowing for arbitrary delaminations.
\newblock {\em International Journal for Numerical Methods in Engineering},
  58(13):2013--2040, 2003.

\bibitem{doi:10.1142/S0219876206001181}
T.~Nagashima and H.~Suemasu.
\newblock Stress analyses of composite laminate with delamination using
  {X-FEM}.
\newblock {\em International Journal of Computational Methods},
  03(04):521--543, 2006.

\bibitem{CurielSosa2012788}
J.L.~Curiel Sosa and N.~Karapurath.
\newblock Delamination modelling of {GLARE} using the extended finite element
  method.
\newblock {\em Composites Science and Technology}, 72(7):788 -- 791, 2012.

\bibitem{Nguyen2013}
V.~P. Nguyen and H.~Nguyen-Xuan.
\newblock High-order {B-splines} based finite elements for delamination
  analysis of laminated composites.
\newblock {\em Composite Structures}, 102:261--275, 2013.

\bibitem{lshape}
B.~G{\"o}zl{\"u}kl{\"u} and D.~Coker.
\newblock Modeling of the dynamic delamination of {L-shaped} unidirectional
  laminated composites.
\newblock {\em Composite Structures}, 94(4):1430 -- 1442, 2012.

\bibitem{Wimmer20082332}
G.~Wimmer and H.E. Pettermann.
\newblock A semi-analytical model for the simulation of delamination in
  laminated composites.
\newblock {\em Composites Science and Technology}, 68(12):2332 -- 2339, 2008.

\bibitem{piegl_book}
L.~A. Piegl and W.~Tiller.
\newblock {\em The {NURBS} Book}.
\newblock Springer, 1996.

\bibitem{rhino}
Rhino.
\newblock {CAD} modeling and design toolkit.
\newblock \url{www.rhino3d.com}.

\bibitem{borden_isogeometric_2011}
M.~J. Borden, M.~A. Scott, J.~A. Evans, and T.~J.~R. Hughes.
\newblock Isogeometric finite element data structures based on {B}\'{e}zier
  extraction of {NURBS}.
\newblock {\em International Journal for Numerical Methods in Engineering},
  87(1‐5):15--47, 2011.

\bibitem{scott_isogeometric_2011}
M.~A. Scott, M.~J. Borden, C.~V. Verhoosel, T.~W. Sederberg, and T.~J.~R.
  Hughes.
\newblock Isogeometric finite element data structures based on {B}\'{e}zier
  extraction of {T}-splines.
\newblock {\em International Journal for Numerical Methods in Engineering},
  88(2):126--156, 2011.

\bibitem{PhuInterface}
V.~P. Nguyen.
\newblock On some practical aspects of fracture modelling using interface
  elements.
\newblock Technical report, Delft University of Technology, The Netherlands,
  2009.
\newblock http://www.academia.edu/1188763/interface-mesh-generator.

\bibitem{Camanho01082003}
P.~P. Camanho, C.~G. Davila, and M.~F. de~Moura.
\newblock Numerical simulation of mixed-mode progressive delamination in
  composite materials.
\newblock {\em Journal of Composite Materials}, 37(16):1415--1438, 2003.

\bibitem{Turon20061072}
A.~Turon, P.P. Camanho, J.~Costa, and C.G. DÃ¡vila.
\newblock A damage model for the simulation of delamination in advanced
  composites under variable-mode loading.
\newblock {\em Mechanics of Materials}, 38(11):1072 -- 1089, 2006.

\bibitem{Benzeggagh1996439}
M.L. Benzeggagh and M.~Kenane.
\newblock Measurement of mixed-mode delamination fracture toughness of
  unidirectional glass/epoxy composites with mixed-mode bending apparatus.
\newblock {\em Composites Science and Technology}, 56(4):439 -- 449, 1996.

\bibitem{vanderMeer2010719}
F.P. van~der Meer and L.J. Sluys.
\newblock Mesh-independent modeling of both distributed and discrete matrix
  cracking in interaction with delamination in composites.
\newblock {\em Engineering Fracture Mechanics}, 77(4):719 -- 735, 2010.

\bibitem{frans}
F.P. van~der Meer.
\newblock Mesolevel modeling of failure in composite laminates: Constitutive,
  kinematic and algorithmic aspects.
\newblock {\em Archives of Computational Methods in Engineering},
  19(3):381--425, 2012.

\bibitem{gutirrez_energy_2004}
M.~A. Guti{\'e}rrez.
\newblock Energy release control for numerical simulations of failure in
  quasi-brittle solids.
\newblock {\em Communications in Numerical Methods in Engineering},
  20(1):19--29, 2004.

\bibitem{verhoosel_dissipation-based_2009}
C.~V. Verhoosel, J.~J.~C. Remmers, and M.~A. Guti{\'e}rrez.
\newblock A dissipation-based arc-length method for robust simulation of
  brittle and ductile failure.
\newblock {\em International Journal for Numerical Methods in Engineering},
  77(9):1290--1321, 2009.

\bibitem{jemjive}
E.~J. Lingen and M.~Stroeven.
\newblock Jem/{Jive}-a {C}++ numerical toolkit for solving partial differential
  equations.
\newblock \url{http://www.habanera.nl/}.

\bibitem{Mi01071998}
Y.~Mi, M.~A. Crisfield, G.~A.~O. Davies, and H.~B. Hellweg.
\newblock Progressive delamination using interface elements.
\newblock {\em Journal of Composite Materials}, 32(14):1246--1272, 1998.

\bibitem{Kress2005458}
G.~Kress, R.~Roos, M.~Barbezat, C.~Dransfeld, and P.~Ermanni.
\newblock Model for interlaminar normal stress in singly curved laminates.
\newblock {\em Composite Structures}, 69(4):458 -- 469, 2005.

\bibitem{Roos2007327}
R.~Roos, G.~Kress, M.~Barbezat, and P.~Ermanni.
\newblock Enhanced model for interlaminar normal stress in singly curved
  laminates.
\newblock {\em Composite Structures}, 80(3):327 -- 333, 2007.

\bibitem{nguyen-offset}
V.~P. Nguyen, P.~Kerfriden, S.P.A. Bordas, and T.~Rabczuk.
\newblock An integrated design-analysis framework for three dimensional
  composite panels.
\newblock {\em Computer Aided Design}, 2013.
\newblock submitted.

\bibitem{bazilevs_isogeometric_2010}
Y.~Bazilevs, {V.M.} Calo, {J.A.} Cottrell, {J.A.} Evans, {T.J.R.} Hughes,
  S.~Lipton, {M.A.} Scott, and {T.W.} Sederberg.
\newblock Isogeometric analysis using {T-splines}.
\newblock {\em Computer Methods in Applied Mechanics and Engineering},
  199(5-8):229--263, 2010.

\end{thebibliography}
